%Paper: alg-geom/9310005
%From: nag@imsc.ernet.in (Subhashis Nag)
%Date: Thu, 7 Oct 93 17:38:06 GMT

%%%%%%%%%%%%%%%%%%%%%%%%%%%%%%%%%%%%%%%%%%%%%%%%%%%%%%%%%%%%%%%%%%%%%%%%
%% MAIN BODY OF TEX PAPER (39 pages) FOLLOWS: %%%%%%%%%%%%%%%%%%%%%%%%%%
%% MAILING SECTIONS 1 to 4 IN ONE FILE AND THEN THE REST %%%%%%%%%%%%%%%
%%%%%%%%%%%%%%%%%%%%%%%%%%%%%%%%%%%%%%%%%%%%%%%%%%%%%%%%%%%%%%%%%%%%%%%%
%%TEICHMULLER THEORY, UNIV PERIOD MAPPING AND THE H-HALF SPACE.(TEX)%%%%
%%%% by Subhashis NAG and Dennis SULLIVAN %%%%%%%%%%%%%%%%%%%%%%%%%%%%%%
%%%%%%%%%%%%%%%%%%%%%%%%%%%%%%%%%%%%%%%%%%%%%%%%%%%%%%%%%%%%%%%%%%%%%%%%

%%%%%%First comes the initial macros file sent from IHES%%%%%%%%%%%%%%%%

\def\today{\ifcase\month\or January\or February\or March\or
April\or May\or June\or July\or August\or September\or
October\or November\or December\fi \space\number\day,
\number\year}

\def\a{\alpha}
\def\b{\beta}

\def\g{\gamma}

\def\l{\lambda}

\def\p{\phi}

\def\t{\theta}

\def\vp{\varphi}

\def\D{\Delta}

\def\Si{\Sigma}

\let\nn=\noindent

\font\tenbb=msbm10
\font\sevenbb=msbm8
\font\fivebb=msbm5
\newfam\bbfam
\textfont\bbfam=\tenbb \scriptfont\bbfam=\sevenbb
\scriptscriptfont\bbfam=\fivebb
\def\bb{\fam\bbfam}

\def\CC{{\bb C}}

\def\RR{{\bb R}}

\def\ZZ{{\bb Z}}

\def\part{\partial}
\def\Inf{\infty}

\def\op{\oplus}
\def\ot{\otimes}

\def\bu{\bullet}

\def\ts{\times}

\def\sbs{\subset}

\def\ra{\rightarrow}

\def\hra{\hookrightarrow}

\def\lbc{\lbrace}
\def\rbc{\rbrace}
\def\lbk{\lbrack}
\def\rbk{\rbrack}

\def\oo{\overline}
\def\ww{\widetilde}

\def\ds{\displaystyle }

\def\Da{{\cal D}}

\def\Ha{{\cal H}}

\def\and{\mathop{\rm and}\nolimits}

\def\all{\mathop{\rm all}\nolimits}

\def\Aut{\mathop{\rm Aut}\nolimits}

\def\Diff{\mathop{\rm Diff}\nolimits}

\def\for{\mathop{\rm for}\nolimits}

\def\grad{\mathop{\rm grad}\nolimits}

\def\Hom{\mathop{\rm Hom}\nolimits}

\def\Hol{\mathop{\rm Hol}\nolimits}

\def\Im{\mathop{\rm Im}\nolimits}

\def\Jac{\mathop{\rm Jac}\nolimits}

\def\Maps{\mathop{\rm Maps}\nolimits}

\def\or{\mathop{\rm or}\nolimits}
\def\of{\mathop{\rm of}\nolimits}

\def\space{\mathop{\rm space}\nolimits}

\def\sgn{\mathop{\rm sgn}\nolimits}

\def\th{\mathop{\rm th}\nolimits}

\catcode`\@=11
\def\displaylinesno #1{\displ@y\halign{
\hbox to\displaywidth{$\@lign\hfil\displaystyle##\hfil$}&
\llap{$##$}\crcr#1\crcr}}

\def\ldisplaylinesno #1{\displ@y\halign{
\hbox to\displaywidth{$\@lign\hfil\displaystyle##\hfil$}&
\kern-\displaywidth\rlap{$##$}
\tabskip\displaywidth\crcr#1\crcr}}
\catcode`\@=12

\def\buildrel#1\over#2{\mathrel{
\mathop{\kern 0pt#2}\limits^{#1}}}

\def\build#1_#2^#3{\mathrel{
\mathop{\kern 0pt#1}\limits_{#2}^{#3}}}

\def\hfl#1#2{\smash{\mathop{\hbox to 6mm{\rightarrowfill}}
\limits^{\scriptstyle#1}_{\scriptstyle#2}}}

\def\hfll#1#2{\smash{\mathop{\hbox to 6mm{\leftarrowfill}}
\limits^{\scriptstyle#1}_{\scriptstyle#2}}}

\def\up#1{\raise 1ex\hbox{\sevenrm#1}}

\def\xx{\vrule height 0.7em depth 0.2em width 0.5 em}

\def\per{|\!\raise -4pt\hbox{$-$}}
 \def\cqfd{\unskip\kern
6pt\penalty 500 \raise
-2pt\hbox{\vrule\vbox to10pt{\hrule
width 4pt \vfill\hrule}\vrule}\par}

\def\trait{\hbox to 12mm{\hrulefill}}
\def\2{{\mathop{\rm
I }\nolimits}\!{\mathop{\rm
I}\nolimits}}

\def\1{{\mathop{\rm I }\nolimits}}

\def\og{\leavevmode\raise.3ex\hbox{\$
scriptscriptstyle\langle\!\langle$}}
\def\fg{\leavevmode\raise.3ex\hbox{\$
scriptscriptstyle \,\rangle\!\rangle$}}

\def\picture #1 by #2 (#3)
{\vcenter{\vskip #2
\special{picture #3}
\hrule width #1 height 0pt depth 0pt
\vfil}}

\def\scaledpicture #1 by #2 (#3 scaled #4){{
\dimen0=#1 \dimen1=#2
\divide\dimen0 by 1000 \multiply\dimen0 by #4
\divide\dimen1 by 1000 \multiply\dimen1 by #4
\picture \dimen0 by \dimen1 (#3 scaled #4)}}

\def\[{{[\mkern-3mu [}}
\def\]{{]\mkern-3mu ]}}

\def\tvi{\vrule height 12pt depth 5pt width 0pt}

\def\cc#1{\hfill\kern .7em#1\kern .7em\hfill}

\def\TeX{T\kern-.1667em\lower.5ex\hbox{E}\kern-.125em X}

\def\ins{{ \raise -2mm\hbox{$<$}
\atop \raise 2mm\hbox{$\sim$}J}
}

\def\sus{{ \raise -2mm\hbox{$>$}
\atop \raise 2mm\hbox{$\sim$}J}
}

%%%%%%%%%%%%%%%%%%%%%%%%%%%%%%%%%%%%%%%%%%%%%%%%%%%%%%%%%%%%%%%%%%%%

\magnification=1200
\overfullrule=0mm
\baselineskip=12.4pt

\vglue 0.4cm
\centerline{\bf{TEICHM\"ULLER THEORY AND THE UNIVERSAL
PERIOD MAPPING}}
\smallskip
\centerline{\bf{VIA QUANTUM CALCULUS AND THE
$H^{1/2}$ SPACE ON THE CIRCLE}}

\medskip
\centerline {by}
\medskip
\centerline{{\bf Subhashis Nag} and {\bf Dennis Sullivan}}

\bigskip

\nn
{\bf Abstract:} Quasisymmetric homeomorphisms
of the circle, that arise in the
Teichm\"uller theory of Riemann surfaces as boundary values of
quasiconfomal diffeomorphisms of the disk, have fractal graphs in
general and are consequently not so amenable to usual
analytical or calculus procedures. In this paper we make use of
the remarkable fact this group $QS(S^{1})$ acts by substitution
(i.e., pre-composition) as a family of
bounded symplectic operators on the Hilbert space $\Ha$=``$H^{1/2}$''
(comprising functions mod constants on $S^1$ possessing
a square-integrable half-order derivative).
Conversely, and that is also important for our work,
quasisymmetric homeomorphisms
are actually {\it characterized} amongst homeomorphisms of
$S^1$ by the property of preserving the space $\Ha$.

Interpreting $\Ha$ via boundary values as the square-integrable first
cohomology of the disk with the cup product symplectic structure, and
complex structure provided by the Hodge star, we obtain a universal form
of the classical period mapping extending the map of [12] [13] from
$Diff(S^{1})/Mobius(S^{1})$ to all of $QS(S^{1})/Mobius(S^{1})$
-- namely to the entire universal Teichm\"uller space, $T(1)$.
The target space for the period map is the universal Siegel
space of period matrices; that is the space of
all the complex structures on $\Ha$ that are compatible with
the canonical symplectic structure.

Using Alain Connes' suggestion of a quantum differential
$d^{Q}_{J}f = [J,f]$ -- commutator of the multiplication operator with
the complex structure operator -- we obtain in lieu of the problematical
classical calculus a quantum calculus for quasisymmetric homeomorphisms.
Namely, one has operators $\{h,L\}$, $d\circ\{h,L\}$, $d\circ\{h,J\}$,
corresponding to the non-linear classical objects
$log (h^\prime)$, ${h'' \over {h'}}dx$, ${1\over 6}Schwarzian(h)dx^{2}$
defined when $h$ is appropriately smooth. Any one of these objects is a
quantum measure of the conformal distortion of $h$ in analogy with the
classical calculus Beltrami coefficient $\mu$ for a
quasiconformal homeomorphism of the disk. Here $L$ is the smoothing
operator on the line (or the circle) with kernel $log \vert x-y \vert$,
$J$ is the Hilbert transform (which is $d \circ L$ or $L \circ d$), and
$\{h,A\}$ means $A$ conjugated by $h$ minus $A$.

The period mapping and the quantum calculus are related in several ways.
For example, $f$ belongs to $\Ha$ if and only if the quantum
differential is Hilbert-Schmidt. Also, the
complex structures $J$ on $\Ha$ lying on the
Schottky locus (image of
the period map) satisfy a quantum integrability condition
$[d^{Q}_{J},J]=0$.

Finally, we discuss the Teichm\"uller space of the universal hyperbolic
lamination ([20]) as a separable complex submanifold
of $T(1)$. The lattice and K\"ahler (Weil-Petersson) metric aspect of the
classical period mapping appear by focusing attention on this space.

\vfill\eject

\baselineskip=14pt
\nn
{\bf \S 1 - Introduction}

The Universal Teichm\"uller Space $T(1)$, which is a
universal parameter space for all Riemann surfaces, is
a complex Banach manifold that may be defined as
the homogeneous space $QS
\left({  S^1}\right)/$M\"ob $\left({  S^1}\right)$.
Here $QS \left({  S^1}\right)$ denotes the group of all
quasisymmetric homeomorphisms of
the unit circle $\left({  S^1}\right)$, and
M\"ob$\left({  S^1}\right)$ is the three-parameter subgroup of M\"obius
transformations of the unit disc (restricted to the
boundary circle). There is
a remarkable homogeneous K\"ahler complex manifold,
$M = \Diff \left({  S^1}\right)/$ M\"ob $\left({ S^1}\right)$,--
arising from applying the Kirillov-Kostant coadjoint orbit method to the
$C^\Inf$-diffeomorphism group
$\Diff \left({  S^1}\right)$ of the circle ([22]) -
that clearly sits embedded in $T(1)$ (since any smooth diffeomorphism
is quasisymmetric).

\bigskip
In [15] it was proved that the canonical complex-analytic
and K\"ahler structures on these two spaces coincide under the natural
injection of $M$ into $T(1)$. (The K\"ahler structure on $T(1)$ is
formal -- the pairing converges on precisely the $H^{3/2}$
vector fields on the circle.) The relevant complex-analytic
and symplectic structures on $M$, (and its close relative $N =
 \Diff \left({  S^1}\right)/ \left({  S^1}\right)$), arise
from the representation theory of
$\Diff \left({  S^1}\right)$ ; whereas on $T(1)$
the complex structure is
dictated by Teichm\"uller theory, and the (formal) K\"ahler
metric is Weil-Petersson.
Thus, the homogeneous space $M$ is a complex analytic
submanifold (not locally closed) in $T(1)$, carrying a canonical
K\"ahler metric.

\bigskip
In subsequent work ([12] [13]) it was shown
that one can canonically associate infinite-dimensional period matrices to
the smooth points $M$ of $T(1)$.  The crucial step in this
construction was a faithful representation (Segal
[18]) of $\Diff \left({S^1}\right)$ on the Frechet space

$$V = C^\Inf \Maps
\left({ S^1, \RR }\right)/\RR
({\mathop{\rm  \  the \ constant \  maps \ }\nolimits})
\eqno (1)$$

\nn
$\Diff \left({  S^1}\right)$  acts by pullback on
the functions in $V$ as a group
of toplinear automorphisms that preserve a basic
symplectic form that $V$ carries.

\bigskip
In order to be able to extend the infinite dimensional
period map to the full space $T(1)$, it is necessary to
replace $V$ by a suitable ``completed'' space
that is invariant under quasisymmetric pullbacks.
Moreover, the quasisymmetric
homeomorphisms should continue to act as bounded
symplectic automorphisms of this extended space. These goals are
achieved in the present paper by developing the theory of
the Sobolev space on the circle consisting of functions with half-order
derivative. This
Hilbert space $ H^{1/2}$ = $\Ha$, {\it which
turns out to be exactly the completion
of the pre-Hilbert space $V$},  actually {\it characterizes}
quasisymmetric (q.s).
homeomorphisms (amongst all homeomorphisms of $S^1$). That fact will be
important for our understanding of the period mapping.
The symplectic structure, $S$, on $V$
extends to $\Ha$ and is preserved by the action of $QS(S^1)$, and
indeed we show that this $S$ is
the {\it unique} symplectic structure available which
is invariant under even the tiny finite-dimensional subgroup M\"ob
$\left({  S^1}\right) $ ($\sbs QS \left({  S^1}\right) $).

\bigskip
We utilise several different characterisations of $\Ha$ and its
complexification. In particular, $\Ha$ comprises
functions on $S^1$ which are
defined quasi-everywhere (i.e., off some set of
logarithmic capacity zero); alternatively, they
appear as non-tangential limits of harmonic functions
of finite Dirichlet energy in the disc.
The last-mentioned fact allows us to interpret $\Ha$ as the
first cohomology space with real coefficients of the unit disc in the
Hodge-theoretic sense. That is important for our subsequent discussion
of the period mapping as a theory of the variation of $S$-compatible
complex structure on this real Hilbert
space $\Ha$. The fact that quasisymmetric
homeomorphisms are the {\it only} ones preserving $\Ha$ is necessary in
our determination of the universal Schottky locus -- namely the image of
$\Pi$.

\bigskip
We present a section
where we discuss quantum calculus on the line (motivated by Alain Connes),
the idea being firstly to demonstrate that
the $H^{1/2}$ functions have such an interpretation. That then allows us
to interpret the universal Siegel space that is the target space for the period
mapping as "almost complex structures on the line" and the Teichm\"uller
points (i.e., the Schottky locus ) can be interpreted as comprising
precisely the subfamily of those complex structures that are {\it
integrable}.

\bigskip
Notice that the fact that capacity zero sets are preserved by
quasisymetric transformations --  whereas merely being
measure zero is not a q.s.-invariant
notion --  goes to exemplify how deeply quasisymmetry
is connected to the properties of $\Ha$. Other characterisations found
below for the complexification $\CC \ot \Ha$ in terms of boundary
values for holomorphic and anti-holomorphic functions are of
independent importance,
and relate to the proof of the uniqueness of the invariant
symplectic structure. That proof utilises a pair of irreducible unitary
representations from the discrete series for
$SL (2, \RR)$ and a version of Schur's lemma.

\bigskip
In universal Teichm\"uller space there resides the separable complex
submanifold $T(H_\Inf)$ -- the Teichm\"uller space of the universal
hyperbolic lamination -- that is exactly the closure of the union of all
the classical Teichm\"uller spaces of closed Riemann surfaces in $T(1)$
(see [20]).
Genus-independent constructions like the universal period mapping proceed
naturally to live on this completed version
of the classical Teichm\"uller spaces. We show that
$T(H_\Inf)$ carries a natural convergent Weil-Petersson pairing.

\bigskip
We make no great claim to originality in this work. Our purpose is to
survey from various different aspects the elegant role of $H^{1/2}$
in universal Teichm\"uller theory, the main goal being to understand
the period mapping in the universal context.
The Hilbert space $\Ha$, its complexification, its symplectic
form and its polarizations etc. appear so naturally in what follows that
it may not be merely facetious to say that the connection of $\Ha$
with Teichm\"uller theory and quasiconformal mappings are not
only ``natural'' but almost ``supernatural''.

\bigskip
\nn
{\bf Acknowledgements:}
It is our pleasant duty to acknowledge gratefully several
stimulating conversations
with Graeme Segal, Michel Zinsmeister, M.S. Narasimhan, Alain Connes,
Ofer Gabber, Stephen Semmes and Tom Wolff.
We heartily thank Michel Zinsmeister for supplying us
with his notes on $H^{1/2}$, and for
generously permitting us to utilise them in this
publication. Also, Graeme Segal suggested the use
of Schur's Lemma in order to
show the essentially unique nature of the symplectic structure.
The ideas on the generalised Jacobi variety (Section 5) arose in
conversations with M.S. Narasimhan.

\medskip
One of us (S.N.) would like to thank very much the
IHES (Bures-sur-Yvette) - where
most of these results were obtained - as well as
the ICTP (Trieste) and the CUNY Graduate Center (New York), for
their excellent hospitality during the Autumn/Winter of 1992.
He would also like to thank the many intersted participants of
the Mathematical Society of Japan International Research Institute on the
``Topology of the moduli space of curves'' (RIMS, Kyoto 1993), for
their discussions and useful feedback.

\bigskip
\nn
{\bf \S - 2 The Hilbert space $H^{1/2} $ on the circle and the line.}

Let $\D$ denote the open unit disc and $U$ the upper half-plane
in the plane ($\CC$), and $S^1 = \part \D$ be the unit circle.

\bigskip
Intuitively speaking, the real Hilbert space
under concern:

$$\Ha \equiv H^{1/2} \left({  S^1, \RR}\right)
 /\RR \eqno (2)$$

\nn
is the subspace of  $L^2 \left({  S^1}\right)$
comprising real functions of mean-value
zero on $S^1$ which have  a half-order derivative also in
$L^2 \left({  S^1}\right)$. Harmonic analysis will tell us that these
functions are actually defined off some set of capacity zero
(i.e., "quasi-everywhere") on the circle, and that they also appear as
the boundary values of real harmonic functions of finite Dirichlet
energy in $\D$. Our first way (of several) to make this precise is
to {\it identify $\Ha$ with the sequence space}

$$\ell_2^{1/2} =
\{
{\mathop{\rm complex \ sequences  }\nolimits}
\ \ u \equiv
(u_1, u_2, u_3, \cdots ): \{
\sqrt{n } \ u_n \}
\
{\mathop{\rm is \ square \ summable \  }\nolimits} \}.
\eqno (3)$$

\bigskip
The identification between (2) and (3) is by showing
 (see Proof of Theorem 2.1)
that the Fourier series

$$f
\left({ e^{i \t } }\right) =
\sum_{ n = - \Inf }^{ \Inf } u_n e^{i n \t} ;
{}~~~u_{- n} = \oo u_n, \eqno (4)$$

\nn
converges quasi-everywhere and defines a real function
of the required type.
The norm on $\Ha$ and on $\ell_2^{1/2}$ is, of course,
 the $\ell_2$ norm of
$\{ \sqrt{n } \ u_n \}$, i.e.,

$$
\left\Vert{ f }\right\Vert_{\Ha}^{2} =
\left\Vert{ u }\right\Vert_{\ell_2^{1/2}}^{2} =
2 \sum_{ n = 1}^{ \Inf  }
n \left\vert{ u_n }\right\vert_{}^{2} .\eqno (5)$$

\bigskip
Naturally $\ell_2^{1/2}$ and $\Ha$ are isometrically
isomorphic separable Hilbert
spaces. Note that $\Ha$ is a subspace of
 $L^2 \left({  S^1}\right)$ because
$\{ \sqrt{n } \ u_n \}$ in $\ell_2$ implies
$\{  \ u_n \} $ itself is in $\ell_2$.

\bigskip
At the very outset let us note the fundamental fact
that the space $\Ha$ is
evidently closed under {\it Hilbert transform }
(``conjugation'' of Fourier series):

$$ (Jf) (e^{i \t}) = -
\sum_{ n = - \Inf }^{ \Inf } i
\sgn (n) u_n e^{i n \t} .\eqno (6)$$

\bigskip

In fact, $J: \Ha \ra \Ha$ is an isometric isomorphism
whose square is the
negative identity, and thus $J$ defines a {\it canonical
 complex structure for $\Ha$}.

\bigskip
\nn
{\bf Remark:}
In the papers [8],[12] [13] [15], we had made use of the fact
that the Hilbert transform defines the
almost-complex structure operator for the tangent space
of the coadjoint orbit manifolds ($M$ and $N$), as well
as for the universal Teichm\"uller space $T(1)$.

\bigskip
Whenever convenient we will pass to a description of our Hilbert space
$\Ha$ as functions on the real line, $\bf R$. This is done by simply
using the M\"obius transformation of the circle onto the line that is
the boundary action of the Riemann mapping ("Cayley transform")
of $\D$ onto $U$. We thus get an isometrically isomorphic copy, called
$H^{1/2}(\RR)$, of our Hilbert space $\Ha$ on the circle
defined by taking $f \in \Ha$ to correspond to $g \in
H^{1/2} (\RR)$ where $g = f \circ R,
R(z) = {z - i \over z + i} $ being the Riemann mapping. The Hilbert
transform complex
structure on $\Ha$ in this version is then described by the usual
singular integral operator on the real line with the
"Cauchy kernel" $(x - y)^{-1}$.

\bigskip
Fundamental for our set up is the dense subspace
 $V$ in $\Ha$ defined by
equation (1) in the introduction. At the level of Fourier
series, $V$ corresponds
to those sequence $\{ u_n \}$ in
$\ell_2^{1/2}$ which go to zero more rapidly than
$n^{- k}$ for any $k >  0$. This
is so because a $C^k$ function has Fourier
ceofficients decaying at least as fast
as $n^{-k}$. Since trigonometric polynomials are in $V$,
it is obvious that $V$ is
norm-dense in $\Ha$. On $V$ one has the basic symplectic
form that we utilised
crucially in [12], [13]:

$$S: V \ts V \ra \RR \eqno (7)$$

\nn
given by

$$S(f, g) =
{1 \over 2 \pi}
\int_{ S^1}^{ } f \cdot dg .
\eqno (8)$$

\nn
This is essentially the signed area of the $(f(e^{i\t}), g(e^{i\t}))$
curve in Euclidean plane.
On Fourier coefficients this bilinear form becomes

$$S(f, g) = 2 \Im\left({ \sum_{n = 1 }^{\Inf} n u_n \oo v_n }\right)
= -i\sum_{n = - \Inf }^{\Inf} n u_n  v_{- n }
\eqno (9)$$

\nn
where $\{ u_n \}$ and $\{ v_n \}$ are respectively
 the Fourier ceofficients of the
(real-valued) functions $f$ and $g$, as in (4).
 Let us note that the
Cauchy-Schwarz inequality applied to (9) shows that
 {\it this non-degenerate
bilinear alternating form extends from $V$ to the full
Hilbert space} $\Ha$.
We will call this extension also
 $S: \Ha \ts \Ha \ra \RR$. Cauchy-Schwarz asserts:

$$
\left\vert{ S(f, g) }\right\vert \leq
\left\Vert{ f }\right\Vert
\cdot
\left\Vert{ g }\right\Vert .
\eqno (10)$$

\nn
Thus $S$ is a jointly continuous - in fact
analytic - map on $\Ha \ts \Ha$.

\bigskip
The important interconnection between the inner
product on $\Ha$, the Hilbert-transform complex
structure $J$, and the form $S$ is encapsulated
in the identity:

$$S \left({ f, Jg }\right) =
\left\langle{  f, g }\right\rangle, \
{\mathop{\rm  for \ all \ }\nolimits}
f, g \in \Ha
\eqno (11)$$

\bigskip
{\it We thus see that $V$ itself was naturally a pre-Hilbert space with
respect to the canonical inner product arising from its symplectic form
and its Hilbert-transform complex structure, and we have just
established that the completion of $V$ is nothing other than the Hilbert
space $\Ha$. Whereas $V$ carried the $C^{\Inf}$ theory, beacause it was
diffeomorphism invariant, the completed Hilbert space $\Ha$ allows
us to carry through our constructions for the full Universal
Teichm\"uller Space because it indeed is quasisymmetrically
invariant.}

\bigskip
It will be important for us to {\it complexify} our
spaces since we need to deal
with isotropic subspaces and polarizations.
Thus we set

$$\CC \ot V \equiv V_\CC = C^\Inf \
{\mathop{\rm Maps }\nolimits}
\left({  S^1, \CC }\right) / \CC
$$

$$
\CC \ot \Ha \equiv \Ha_\CC = H^{1/2}
\left({  S^1, \CC }\right) / \CC
\eqno (12)
$$

\nn
$\Ha_\CC$ is a complex Hilbert space isomorphic to
$\ell_2^{1/2} (\CC)$ - the latter comprising
the Fourier series

$$f
\left({ e^{i \t } }\right) =
\sum_{n = - \Inf }^{ \Inf }
u_n e^{i n \t} \ \ , \ \ u_0 = 0
\eqno (13)$$

\nn
with $\{ \sqrt{ |n| } \ u_n \}$ being square summable
 over $\ZZ - \{ 0 \}$. Note
that the Hermitian inner product on $\Ha_\CC $
derived from (5) is given by

$$\left\langle{  f, g}\right\rangle =
\sum_{ n = - \Inf}^{ \Inf  } |n| u_n \oo v_n.
\eqno (14)$$

\nn
[This explains why we introduced the factor 2 in the formula
(5).] The fundamental {\it orthogonal decomposition} of
$\Ha_\CC$ is given by

$$\Ha_\CC = W_+ \op W_- \eqno (15)$$

\nn
where

$$W_+ =
\{ f \in \Ha_\CC: \
{\mathop{\rm  all \ negative \ index \ Fourier \
coefficients \ vanish \
}\nolimits} \}$$

\nn
and
$$\oo W_+ = W_- =
\{ f \in \Ha_\CC: \
{\mathop{\rm  all \ positive \ index \ Fourier \
coefficients \ vanish \
}\nolimits}
\}.$$

\bigskip
Here we denote by bar the complex anti-linear
automorphism of $\Ha_\CC$ given by
conjugation of complex scalars.

\bigskip
Let us extend $\CC$-linearly the form $S$
and the operator $J$ to $\Ha_\CC$ (and
consequently also to  $V_\CC$).
The complexified $S$ is still given by the
right-most formula in (9). Notice that
$W_+ $ and $W_-$ can be characterized as
precisely the $- i$
{\it and $+ i$ eigenspaces} (respectively)
 of the $\CC$-linear extension of $J$, the Hilbert transform.
Further, each of $W_+ $ and $W_-$ is
{\it isotropic} for $S$, i.e.,
$S(f, g) = 0$, whenever both $f$ and $g$
 are from either $W_+$ or $W_-$ (see
formula (9)). Moreover, $W_+$ and $W_-$
are {\it positive isotropic} subspaces in the
sense that the following identities hold:

$$
\left\langle{ f_+, g_+ }\right\rangle_{}^{} =
i S
\left({ f_+, \oo g_+ }\right), \ for \
f_+, g_+ \in W_+ \eqno (16)$$

\nn
and

$$
\left\langle{ f_-, g_- }\right\rangle_{}^{} = -
i S
\left({ f_-, \oo g_- }\right), \ for \
f_-, g_- \in W_- . \eqno (17)$$

\bigskip
\nn
{\bf Remark:}
(16) and (17) show that we could have defined the
 inner product and norm on
$\Ha_\CC$ from the symplectic form $S$, by using these relations to {\it
define } the inner products on $W_+$ and $W_-$, and declaring $W_+$ to be
perpendicular to $W_-$. Thus, for general $f, g \in \Ha_\CC$
one has the fundamental identity

$$
\left\langle{ f, g }\right\rangle_{}^{} =
i S
\left({ f_+, \oo g_+ }\right) -
i S
\left({ f_-, \oo g_- }\right) .
\eqno (18)$$

\bigskip
\nn
We have thus described the Hilbert space
structure of $\Ha$ simply in terms of the
canonical symplectic form it carries and
the fundamental decomposition (15).
[Here, and henceforth, we will let $f_{\pm}$
denote the projection of $f$ to
$W_{\pm}$, etc.].

\bigskip
In order to prove the first results of this paper,
we have to rely on interpreting the functions in $H^{1/2}$
as boundary values (``traces'') of
functions in the disc $\D$ that have finite
Dirichlet energy, (i.e. the first
derivatives are in $L^2 (\D)$). We start
explaining this material.

\bigskip
Define the following ``Dirichlet space'' of
 harmonic functions:

$$\Da =
\{ F : \D \ra \RR : F \
{\mathop{\rm  is \ harmonic \ }\nolimits},
F(0) = 0, \ \and \ E(F) < \Inf \}
\eqno (19)$$

\nn
where the energy $E$ of any (complex-valued) $C^1$ map on $\D$
is defined as the $L^2 (\D) $ norm of $\grad (F)$ :

 $$
\left\Vert{ F }\right\Vert_{\Da}^{2} = E(F) =
{1\over 2 \pi }
\int_{ }^{ }\!
\int_{ \D}^{ }
\left({ \left\vert{ { \part F \over
\part x }  }\right\vert_{}^{2} +
\left\vert{ { \part F \over
 \part y }  }\right\vert_{}^{2}
}\right) dx dy
\eqno (20)$$

\nn
$\Da$, and its complexification $\Da_\CC$,
 will be Hilbert spaces with respect to
this energy norm.

\bigskip
We want to identify the space $\Da$ as precisely
the space of harmonic functions in
$\D$ solving the Dirichlet problem for functions
 in $\Ha$. Indeed, the {\it
Poisson integral representation allows us
to map $P: \Ha \ra \Da$ so that $P$ is
an isometric isomorphism of Hilbert spaces.}

\bigskip
To see this let
$f (e^{i \t}) = \ds
\sum_{ - \Inf }^{ \Inf } u_n e^{i n \t}$
be an arbitrary member of
$\Ha_\CC$. Then the Dirichlet extension of
$f$ into the disc is:

$$ F(z) =
\sum_{n =  - \Inf }^{ \Inf }
u_n r^{|n|} e^{i n \t } =
\left({ \sum_{n = 1 }^{ \Inf } u_n z^{n} }\right) +
\left({ \sum_{m = 1 }^{ \Inf } u_{- m } \oo z^{m}  }\right)
\eqno (21)$$

\nn
where $z = re^{i \t}$.
 From the above series one can directly compute the $L^2 (\D)$
norms of $F$ and also of $\grad (F) = (\part F/ \part x, \part F/ \part y)$.
One obtains the following:

$$E(F) =
{1 \over 2 \pi }
\int_{ }^{ }\!
\int_{ \D}^{ }
| \grad (F) |^2 =
\sum_{ - \Inf }^{ \Inf } | n | |u_n |^2 \equiv
\left\Vert{ f }\right\Vert_{\Ha}^{2}
< \Inf
\eqno (22)$$

We will require crucially the well-known formula of Douglas (see [2,pg.
36-38])
expressing the above energy of $F$ as the double integral on
$S^1$ of the square of the first differences of the boundary values $f$.

$$
E(F)=
{1\over {16{\pi}^2}}
\int_{S^1}^{}\!
\int_{S^1}^{}
[{(f(e^{i\t})-f(e^{i\p}))}/{sin({(\t - \p)}/2)}]^2 d\t d\p
\eqno(23)
$$

Transferring to the real line by the M\"obius transform identification
of $\Ha$ with $H^{1/2}(\RR)$ as explained before, the above identity
becomes as simple as:

$$
E(F)=
{1\over {4{\pi}^2}}
\int_{\RR}^{}\!
\int_{\RR}^{}
[{(f(x)-f(y))}/{(x-y)}]^2 dxdy=
\left\Vert{ f }\right\Vert^{2}
\eqno(24)
$$

\bigskip
\nn
Calculating from the series (21), the $L^2$-norm of $F$ itself is:

$$
{1 \over 2 \pi }
\int_{ }^{ }\!
\int_{\D }^{ }
| F|^2 dx dy
= \sum_{ - \Inf }^{ \Inf }
{  |u_n |^2  \over (|n| + 1)}
\leq E(F) < \Inf
\eqno (25)$$

\nn
(22) shows that indeed
{\it Dirichlet extension is isometric from }
$\Ha $ to $\Da$, whereas (25) shows
that the functions in $\Da$ are themselves in
$L^2$, so that the {\it the inclusion
of  $\Da \hra L^2 (\D)$ is continuous}.
(Bounding the $L^2$ norm of $F$ by the
$L^2$ norm of its derivatives is  a ''Poincar\'e
inequality'').

\bigskip
It is therefore clear that $\Da$
is a subspace of the usual Sobolev space $H^1
(\D)$ comprising those functions
in $L^2 (\D)$ whose first derivatness (in the
sense of distributions) are also in $L^2 ( \D)$.
The theory of function spaces implies
(by the ``trace theorems'') that $H^1$ functions
lose half a derivative in going to
a boundary hyperplane. Thus it is known that the
functions in $\Da$ will indeed have
boundary values in $H^{1/2}$.
See [5] [7] and [21].

\bigskip
Moreover, the identity (24) shows that
for any $F \in \Da$, the Fourier expansion of
the trace on the boundary circle is
 a Fourier series with
$\sum |n| |u_n|^2 < \Inf$.
But Fourier expansions with coefficients in such a
weigted $\ell_2$ space, as in our situation,
 are known to converge {\it
quasi-everywhere} (i.e. off a set of logarithmic
 capacity zero) on $S^1$. See
Zygmund [23, Vol 2, Chap. XIII]. The identification
between $\Da$ and
$\Ha$ (or $\Da_\CC$ and $\Ha_\CC$) is now complete.

\bigskip
It will be necessary for us to identify the
 $W_\pm$ polarization of $\Ha_\CC$ at
the level $\Da_\CC$. In fact, let us decompose
the harmonic function $F$ of (21)
into its holomorphic and anti-holomorphic parts ;
these are $F_+ $ and $F_-$,
which are (respectively) the two sums bracketed
 separately on the right hand side
of (21). Clearly $F_+$ is a holomorphic
function extending $f_+$ (the $W_+$ part
of $f$), and $F_-$ is anti-holomorphic extending
 $f_-$. We are thus led to
introduce the following space of holomorphic functions
whose derivatives are in $L^2 (\D)$:

$$
\Hol_2 (\D) =
\{ H: \D \ra \CC: H \
{\mathop{\rm is \ holomorphic \  }\nolimits},
H(0) = 0
\ \
\and
\ \
\int_{ }^{ }\!
\int_{\D }^{ }
| H' (z) |^2 dx dy < \Inf \} .\eqno (26)$$

\bigskip
\nn
This is a complex Hilbert space with the norm

$$
\left\Vert{ H }\right\Vert_{}^{2}
 = {1 \over 2 \pi }
\int_{ }^{ }\!
\int_{ \D }^{ }
| H' (z)  |^2 dx dy .
\eqno (27)$$

\bigskip
\nn
If $H(z) = \ds \sum_{ n = 1}^{\Inf } u_n z^n$,
 a computation in polar coordinates
(as for (21), (25)) produces

$$
\left\Vert{ H }\right\Vert_{}^{2} =
\sum_{ n = 1}^{\Inf }
n \left\vert{ u_n }\right\vert_{}^{2}.
\eqno (28)$$

\nn
Equations (27) and (28) show that the norm-squared
is the Euclidean area of the
(possibly multi-sheeted) image of the map $H$.

\bigskip
We let $\oo {\Hol_2} (\D)$ denote the Hilbert space
of antiholomorphic functions
conjugate to those in $\Hol_2 (\D)$. The norm is defined by
stipulating that the anti-linear isomorphism of $\Hol_2$ on
$\oo {\Hol_2} $ given by conjugation should be an isometry. The
Cauchy-Riemann equation for $F_+$ and $\oo F_-$ imply that

$$
\left\vert{ \grad (F) }\right\vert_{}^{2}
= 2 \left\lbc{ \left\vert{ F'_+ }\right\vert_{}^{2}
- \left\vert{ F'_- }\right\vert_{}^{2} }\right\rbc
\eqno (29)$$

\nn
and hence
$$
\left\Vert{ F_+ }\right\Vert_{}^{2}
+ \left\Vert{ F_- }\right\Vert_{}^{2} =
\left\Vert{ f }\right\Vert_{\Ha_\CC}^{2}.
\eqno (30)$$

\bigskip
Now, the relation between $\Da$ (harmonic functions in
$H^1 (\D))$ and $\Hol_2 (\D)$ is transparent, so that the
holomorphic functions in $\Hol_2$ will have non-tangential
limits quasi-everywhere on $S^1$ - defining a function $W_+$.
We thus collect together, for the record,
the various representations of our basic Hilbert space:

\bigskip
\noindent
{\bf THEOREM 2.1:}
There are {\it canonical isometric isomorphisms} between the
following complex Hilbert spaces:

\medskip
(1) $\Ha_\CC = H^{1/2} \left({ S^1, \CC }\right) / \CC
= \CC \ot H^{1/2}(\RR)
= W_+ \op W_-$;

\medskip
(2) The sequence space $\ell_2^{1/2} (\CC)$ (constituting the Fourier
coefficients of the above quasi-everywhere defined functions);

\medskip
(3) $\Da_\CC$, comprising normalized finite-energy harmonic
functions (either on $\D$ or on the half-plane $U$); [the norm-squared
being given by (20) or (22) or (23) or (24)];

\medskip
(4) $\Hol_2 (\D) \op \oo {\Hol_2} (\D)$.

\bigskip
\nn
Under the canonical identifications, $W_+$ maps to
$\Hol_2 (\D)$ and $W_-$ onto $ \oo {\Hol_2} (\D)$.
\xx

\bigskip
\noindent
{\bf Remark:}
One advantage of introducing the full Sobolev space $H^1 (\D)$
(rather than only its harmonic subspace $\Da$) is that we may
use {\it Dirichlet's principle} to rewrite the norm on $\Ha$
as

$$
\left\Vert{ f }\right\Vert_{\Ha}^{2} =
\inf
\{ E(F): F \
{\mathop{\rm ranges \ over \ all \
extensions \ to \ }\nolimits} \D \ \of \ f \}
\eqno (31)$$

\bigskip
\noindent
By Dirichlet principle, the infimum is realized by the
harmonic extension $P(f) = F$ of (23). In connection with this
it is worth pointing out still another formula for the norm :

$$
\left\Vert{ f }\right\Vert_{\Ha}^{2} =
\int_{ S^1 }^{ }
F \cdot { \part F \over \part n} ds
\eqno (32)$$

\nn
where $F$ is the harmonic extension to $\D$ of $f$. This
follows from the well-known Green's identity:

$$
\int_{ }^{ }\!
\int_{ \D}^{ }
\left\vert{ \grad F }\right\vert_{}^{2} =
\int_{\D}^{ } \!
\int_{  }^{ }
F (\D F) +
\int_{ S^1 }^{ }
F \cdot
{\part F \over \part n} \cdot ds
\eqno (33)$$

\bigskip
\noindent
The first term on the right drops out since $F$ is harmonic.
Hence (32) follows. The close relation of formula (32) with
the symplectic pairing formula (8) should be noted.

\bigskip
\noindent
{\bf Remark:}
Since the isomorphisms of the Theorem are all isometric, and
because the norm arose from the canonical symplectic structure,
(formulas (16), (17), (18)), it is instructive to work out the
formulas for the symplectic form $S$ on $\Da_\CC$ and on
$\Hol_2 (\D) $. This is left to the interested reader.

\nn
{\bf \S 3 - Quasisymmetric invariance.}

Quasiconformal (q.c.) self-homeomorphisms of the disc $\D$ (of
the upper half-plane) $U$ are known to extend continuously to
the boundary. The action on the boundary circle (respectively,
on the real line $\RR$) is called a {\it quasisymmetric}
(q.s.) homeomorphism. By [4], $\vp: \RR
\ra \RR$ is quasisymmetric if and only if, for all $x \in
\RR$ and all $t > 0$, there exists some $K > 0$ such that

$$
{1 \over K}
\leq { \vp (x + t) - \vp (x) \over  \vp (x) - \vp (x - t) }
\leq K \eqno (34)$$

\bigskip
\noindent
On the circle this condition for $\vp: S^1 \ra S^1$ means
that $|\vp(2I)|/|\vp(I)| \leq K$, where $I$ is any interval on $S^1$
of length less that $\pi$, $2I$ denotes the interval obtained
by doubling $I$ keeping the same mid-point, and $| \bu |$
denotes Lebesgue measure on $S^1$. See  [3] [11] and [14]
as general references.

\bigskip
Given any orientation preserving homeomorphism $\vp: S^1 \ra
S^1$, we use it to pullback functions in $\Ha$ by
precomposition:

$$V_\vp (f) = \vp^* (f) =
f \circ \vp -
{1 \over 2 \pi}
\int_{S^1 }^{ }
(f \circ \vp).
\eqno (35)$$

\bigskip
\noindent
[We subtract off the mean value in order that the resulting
function also possess zero mean. Since constants pullback to
themselves, the operation is well-defined on $\Ha$.]

\bigskip
We prove:
\bigskip
\noindent
{\bf THEOREM 3.1:}
$V_\vp$ maps $\Ha$ to itself (i.e., the space $\Ha \circ
\vp$ is $\Ha$) if and only if $\vp$ is quasisymmetric.
The operator norm of $V_\vp \leq {\sqrt{K + K^{-1}}}$, whenever $\vp$
allows a K-quasiconformal extension into the disc.

\bigskip
\noindent
{\bf COROLLARY 3.2:}
The group of all quasisymmetric homeomorphism on $S^1, QS
\left({ S^1 }\right)$, acts faithfully by bounded toplinear
automorphisms on the Hilbert space $\Ha$ (and therefore also
on $\Ha_\CC$).

\bigskip
\noindent
{\bf Proof of sufficiency}

Assume $\vp$ is q.s. on $S^1$, and let $\Phi:\D \ra \D$ be
any quasiconformal extension. Let $f \in \Ha$ and suppose
$P(f) = F \in \Da$ is its unique harmonic extension into $\D$.
Clearly $ G \equiv F \circ \Phi$ has boundary values $f \circ
\vp$, the latter being (like $f$) also a continuous function
on $S^1$ defined quasi-everywhere. [Here we recall that q.s.
homeomorphisms carry capacity zero sets to again such sets,
although measure zero sets can become positive measure.]
To prove that $f \circ \vp$ minus its mean value is in $\Ha$,
it is enough to prove that the Poisson integral of $f \circ
\vp$ again has finite Dirichelet energy. Indeed we will show

$$E
({\mathop{\rm  harmonic \ extension \ of }\nolimits}
\ \vp^* (f)) \leq 2
\left({ {1 + k^2 \over 1 - k^2} }\right) E (F). \eqno (36)$$

\noindent
Here $0 \leq k < 1$ is the q.c. constant for $\Phi$, i.e.,

$$\left\vert{ \Phi_{\oo z} }\right\vert_{}^{}
\leq k
\left\vert{  \Phi_{ z} }\right\vert_{}^{}
\ \ {\mathop{\rm
a.e. \ in  }\nolimits}  \ \D .
\eqno (37)$$

\noindent
The operator norm of $V_\vp$ is thus no more that
$2^{1/2 }
\left({ {\ds 1 + k^2 \over \ds 1 - k^2} }\right)^{1/2 }$. The last
expression is equal to the bound quoted in the Theorem, where, as usual,
$ K =  (1+k)/(1-k)$.

\bigskip
Towards establishing (36) we prove that the inequality holds
with the left side being the energy of the map $G = F \circ
\Phi$. Since $G$ is therefore also a finite energy extension
of $f \circ \vp $ to $\D$, Dirichlet's principle (see (31)
above) implies the required inequality. (Note that since
Dirichlet integral is insensitive to adding a constant to a
function, the energy of $G$ is the same as the energy of
$G-G(0)$.)

To compute $E(G)$ we note that

$$
\left({ { \part G \over \part x } }\right)^2
+
\left({ { \part G \over \part y } }\right)^2
\leq 2
\left\lbk{
\left({ { \part F \over \part u } }\right)^2
+
\left({ { \part F \over \part v } }\right)^2}\right\rbk
\left\lbk{
\left\vert{ \Phi_{ z} }\right\vert_{}^{2}
+
\left\vert{ \Phi_{\oo z} }\right\vert_{}^{2} }\right\rbk.
\eqno (38)$$

\noindent
We have writen $\Phi(x, y) = u(x, y) + iv(x, y)$ and
$F = F(u, v)$. (38) follows by straight computation using the
chain rule. But notice that the Jacobian of $\Phi$ is

$$\Jac (\Phi) =
\left\vert{ \Phi_{ z} }\right\vert_{}^{2}
-
\left\vert{ \Phi_{\oo z} }\right\vert_{}^{2} .
\eqno (39)$$

\noindent
By the quasiconformality (37) we therefore get from (38):

$$
\left\lbk{ G_x^2 + G_y^2 }\right\rbk
\leq 2
\left({ { 1 + k^2 \over  1 - k^2} }\right)
\left\lbk{ F_u^2 + F_v^2 }\right\rbk
\Jac (\Phi)  \eqno (40)$$

\nn
Using change of variables in the Dirichlet integral we
therefore derive

$$E(G) \leq 2
\left({ { 1 + k^2 \over  1 - k^2} }\right) E (F)
\eqno (41)$$

\nn
as desired. \xx

\bigskip
\noindent
{\bf Remark:}
Since the Dirichlet integral in two dimensions is
invariant under conformal mappings, it is
not too surprising that it is quasi-invariant
under quasiconformal transformations. Such quasi-invariance
has been noted before and is applied, for example, in
[1] and [16].

\bigskip
\noindent
{\bf Proof of necessity:}
As we mentioned before, the idea of this proof is
taken from the notes of M. Zinsmeister.

\bigskip
Since two-dimensional Dirichlet integrals are conformally
invariant, we will pass to the upper half-plane $U$ and its
boundary line $\RR$ to aid our presentation. As explained earlier, using
the Cayley transform we transfer
everything over to the half-plane; the traces
on the boundary constitute the space of quasi-everywhere defined
functions called  $H^{1/2} (\RR)$.

 From the Douglas identity, equation (24), we recall
that an {\it equivalent way of expressing the Hilbert space norm on }
$H^{1/2} (\RR)$ is

$$
\left\Vert{ g }\right\Vert_{}^{2}
=
{1 \over {4 {\pi^2}}}
\int_{ }^{ }
\!\int_{\RR^2 }^{ }
\left\lbk{
{ g(x) - g(y) \over x - y }  }\right\rbk_{}^{2}
dx dy,
{}~~g \in H^{1/2} (\RR). \eqno (42)$$

\noindent
Equation (42) immediately shows that
 $ \left\Vert{ g }\right\Vert = \left\Vert{ \ww g }\right\Vert$
where $\ww g(x) = g(ax + b)$ for any real $a (\ne 0) $ and
$b$. This will be important.

Assume that $\vp: \RR \ra \RR$ is an orientation preserving
homeomorphism such that
$V_{\vp^{-1}}: H^{1/2} (\RR) \ra H^{1/2} (\RR)$ is a bounded
automorphism.
Let us say that the norm of this operator is $M$.

Fix a bump function $f \in C_0^\Inf (\RR)$ such that $f \equiv
1$ on $[- 1, 1] f \equiv 0$ outside $[- 2, 2]$ and $0 \leq f
\leq 1$ everywhere. Choose any $c \in \RR$ and any positive
$t$. Denote $I_1 = [x - t, x]$ and $I_2 = [x, x + t]$. Set
$g(u) = f(au + b)$, choosing $a$ and $b$ so that $g$ is
identically $1$ on $I_1$ and zero on $[x + t, \Inf)$.

By assumption, $g \circ \vp^{- 1}$ is in $H^{1/2} (\RR)$ and
$
\left\Vert{ g \circ \vp^{-1} }\right\Vert_{}^{}
\leq M \left\Vert{ g }\right\Vert_{}^{} =
M \left\Vert{ f }\right\Vert_{}^{}$.
We have

$$
\eqalignno{
 M \left\Vert{ f }\right\Vert_{}^{}
&\geq
\int_{ }^{ }
\!\int_{\RR^2 }^{ }
\left\lbk{
{ g \circ \vp^{- 1}(u) -
g \circ \vp^{- 1}(v) \over u - v }  }\right\rbk_{}^{2}
du dv\cr
&\geq
\int_{ v = \vp (x - t) }^{ v = \vp (x)}
\int_{ u = \vp (x + t) }^{ u = \Inf}
{1 \over (u - v)^2 } du dv\cr
&=
\log
\left({ 1 +  {  \vp (x) - \vp (x - t) \over
\vp (x + t) - \vp (x)  } }\right)  .
&(43) \cr}$$

\nn
[We utilise the elementary integration
$\ds \int_{\g}^{ \Inf}
\ds \int_{ \a }^{ \b}
{1 \over (u - v)^2 } du dv =
\log
\left({ 1 + {\b - \a \over \g - \b } }\right) $, for
$\a<\b<\g$. ] We thus obtain the result that

$${  \vp (x + t ) - \vp (x ) \over
\vp (x ) - \vp (x - t)  }
\geq
{ 1 \over e^{M \left\Vert{ f }\right\Vert } - 1 }
\eqno (44)$$

\nn
for arbitrary real $x$ and positive $t$. By utilising
symmetry, namely by shifting the bump to be $1$ over $I_2$ and
$0$ for $u \leq x - t$, we get the opposite inequality:

$$
{  \vp (x + t ) - \vp (x ) \over
\vp (x ) - \vp (x - t)  }
\leq  e^{M \left\Vert{ f }\right\Vert } - 1.
\eqno (45)$$

\bigskip
\nn
The Beurling-Ahlfors condition on $\vp$ is verified, and we
are through. Both the theorem and its corollary are proved. \xx

\bigskip
\nn
{\bf \S 4 - The invariant symplectic structure.}

The quasisymmetric homeomorphism group,
$QS \left({  S^1}\right)$, acts on $\Ha$ by precomposition
(equation (35)) as bounded operators, {\it preserving the
canonical symplectic form}
$S: \Ha \ts \Ha \ra \RR$ (introduced in (8), (9), (10)).
This is the central fact which we will analyse in this
section. It is the crux on which the extension of the period
mapping to all of $T(1)$ hinges:

\bigskip
\nn
{\bf PROPOSITION 4.1:}
For every $\vp \in QS \left({  S^1}\right)$, and all $f, g \in
\Ha$,

$$S
\left({ \vp^* (f), \vp^* (g)  }\right) =
S (f, g). \eqno (46)$$

\nn
Considering the complex linear extension of the action to
$\Ha_\CC$, one can assert that the only quasisymmetrics which
preserve the subspace $W_+ = \Hol_2 (\D)$ are the M\"obius
transformations. Then M\"ob $\left({  S^1}\right)$ acts as
unitary operators on
$W_+$ (and $W_-$).

\bigskip
Before proving the proposition we would like to point out that
this canonical symplectic form enjoys a far stronger
invariance property:

\bigskip
\noindent
{\bf LEMMA 4.2:}
If $\vp: S^1 \ra S^1$ is any (say $C^1$) map of winding
number (= degree) $k$, then

$$S(f \circ \vp, g \circ \vp) = kS (f, g) \eqno (47)$$

\nn
for arbitrary choice of ($C^1$) functions $f$ and $g$ on
the circle. In particular, $S$ is invariant under pullback by
all degree one mappings.

\bigskip
\noindent
{\bf Proof:}
The proof of (47), starting from (8), is an exercise in
calculus. Lift $\vp$ to the universal cover to obtain $\ww \vp
: \RR \ra \RR $; the degree of $\vp$ being $k (\in \ZZ)$
implies that
$\ww \vp (t + 2 \pi) = \ww \vp (t) + 2 k \pi$. Breaking up
$[0, 2 \pi]$ into pieces on which $\ww \vp$ is monotone, and
applying the change of variables formula in each piece,
produces the result. \xx

\bigskip
\noindent
{\bf Proof of Proposition 4.1:}
The Lemma shows that (46) is true whenever the quasisymmetric
homeomorphism $\vp$ is at least $C^1$. By
Lehto-Virtanen [11, Chapter II, Section 7.4] we know that for
arbitrary q.s $\vp$, there exist real analytic q.s.
homeomorphisms $\vp_m$ (with the same quasisymmetry constant
as $\vp$) that converge uniformly to $\vp$. An approximation
argument, as below, then proves the required result.

Let us denote the $n^{\th}$ Fourier coefficient of a function
$f$ on $S^1$ by $F_n (f)$. Recall from equation (9) that

$$S(f, g) = - i \sum_{ n = - \Inf }^{ \Inf } nF_n (f) F_{-n} (g)
\eqno (48)$$

\nn
for all $f, g$ in $\Ha_\CC$. Now since $S$ is continuous it is
enough to check (46) on the dense subspace $V$ of smooth
functions $f$ and $g$. Therefore assume $f$ and $g$ to be
smooth.

Since $\vp_m \ra \vp $ uniformy it follows that $F_n (f \circ
\vp_m) \ra F_n (f \circ \vp) $ as $ m \ra \Inf$ (for each fixed
$n$). Applying the dominated convergence theorem to the sums
(48) we immediately see that as $m \ra \Inf$,

$$ S
\left({ \vp_m^* (f), \vp_m^* (g) }\right)
\ra
 S
\left({ \vp^* (f), \vp^* (g) }\right).
\eqno (49)$$

\nn
Lemma 4.1 says that for each $m$,
$S\left({ \vp_m^* (f), \vp_m^* (g) }\right) = S(f, g)$.
We are through.

\bigskip
If the action of $\vp$ on $\Ha_\CC$ preserves $W_+$ it is easy
to see that $\vp$ must be the boundary values of some
holomorphic map $\Phi: \D \ra \D$. Since $\vp$ is a
homeomorphism one can see that $\Phi$ is  a holomorphic
homeomorphism (as explained also in [12, Lemma of Section 1])
- hence a M\"obius transformation. Since every $\vp$ preserves
$S$, and since $S$ induces the inner product on $W_+$ and
$W_-$ by (16) (17), we note that such a symplectic
transformation preserving $W_+$ must necessarily act unitarily.
\xx

\bigskip
\noindent
{\bf Remark:}
The remarkable invariance property (47) leads us to ask a
question that may shed light on the structure of degree $k$
maps of $S^1$ onto itself.
Given a vector space $V$ equipped with a bilinear form $S$,
one may fix some constant $k (\ne 0)$ and study the family of
linear maps $A$ in $\Hom (V, V) $ such that

$$S\left({ A \left({  v_1}\right) ,
A \left({  v_2}\right)
}\right) =k S\left({  v_1, v_2 }\right)
\eqno (50)$$

\nn
holds for all  $v_1, v_2$ in $V$. Of course, the trivial
multiplication (by $ \sqrt{ k} $) will be such a map, but we
have in Lemma 4.1 a situation where the interesting family of
linear maps obtained by degree $k$ pullbacks provide a profusion
of examples - precisely when $k$ is an integer.

Furthermore, in the situation at hand, we may take $V$ as the
space of $C^\Inf$ (real or complex) functions on the circle. Then $V$
also carries {\it algebra} structure by pointwise multiplication.
The pullbacks by degree $k$ mappings clearly preserve
this multiplicative structure
(whereas dilatations do not). It is interesting to question
whether the only linear maps that preserve the algebra
structure and also satisfy the relation (50), (for integer
$k$), must necessarily arise from some degree $k$ mapping of
$S^1$ on itself.

\bigskip
Theorem 3.1 and Proposition 4.1 enable us to consider $QS
\left({  S^1}\right)$ as a subgroup of the bounded symplectic
operators on $\Ha$. Since the heart of the matter in extending
the period mapping from Witten's homogeneous space $M$ (as in
[12], [13]) to $T(1)$ lies in the property of preserving this
symplectic form on $\Ha$, we prove below that $S$ is indeed the
{\it unique} symplectic form that is $\Diff
\left({  S^1}\right)$ or $QS \left({  S^1}\right)$ invariant.
It is all the more surprising that the form $S$ is canonically
specified by requiring its invariance under simply the
3-parameter subgroup M\"ob
$\left({  S^1}\right) (\hra \Diff \left({  S^1}\right)
\hra QS\left({  S^1}\right)$ ).

\bigskip
\noindent
{\bf THEOREM 4.3:}
Let $S \equiv S_1$ be the canonical symplectic form on $\Ha$.
Suppose $S_2: \Ha \ts \Ha \ra \RR$ is any other continuous
bilinear form such that
$S_2 ( \vp^* (f), \vp^* (g)) = S_2 (f, g)$, for all $f, g$ in
$\Ha$ whenever $\vp$ is in M\"ob
$\left({  S^1}\right) $. Then $S_2$ is necessarily a real
multiple of $S$. Thus every form on $\Ha$ that is M\"ob
$\left({  S^1}\right) \equiv PSL (2, \RR)$ invariant is
necessarily non-degenerate (if not identically zero) and
remains invariant under the action of the whole of $QS
\left({  S^1}\right)$. (Also, it automatically satisfies the
even stronger invariance property (47)).

\bigskip
The proof requires some representation theory. Since this
paper is written with complex analysts in mind, we have presented
some detail. We start with:

\bigskip
\noindent
{\bf LEMMA 4.4:}
The duality induced by canonical form $S_1$ is (the negative
of) the Hilbert transform (equation (6)). Thus the map $\Si_1$
(induced by $S_1$) from $\Ha$ to $\Ha^*$ is an invertible
isomorphism.

\noindent
{\bf Proof:}
Given the continuous bilinear pairing $S_i: \Ha \ts \Ha \ra
\RR (i = 1, 2)$ we are considering the induced ``duality'' maps

$$
\Si_i: \Ha \ra \Ha^* \ \ \ \  (i = 1, 2) \eqno (51)$$

\nn
which are bounded linear operators defined by:

$$
\Si_i (g) = S_i (\bu, g) \ , \ g \in \Ha. \eqno (52)$$

\nn
since $\Ha$ is a Hilbert space, the dual $\Ha^*$ is canonically
isomorphic to $\Ha$ via:

$$\l(g) =
\left\langle{\bu, g }\right\rangle.
\eqno (53)$$

\nn
Here $\l(g)$ is the linear functional represented by $g \in
\Ha$. Equation (11) says
$S_1 (f, Jg) =
\left\langle{ f, g }\right\rangle$ and therefore that

$$\Si_1 (g) = \l (- Jg) \eqno (54)$$

\nn
as required. \xx

\bigskip
The basic tool in proving the Theorem 4.3 is to consider the
{\it ``interwining operator''}

$$M = \Si_1^{-1} \circ \Si_2: \Ha \ra \Ha \eqno (55)$$

\nn
which is a bounded linear operator on $\Ha$ by the above Lemma.

\bigskip
\noindent
{\bf LEMMA 4.5:}
$M$ commutes with every invertible linear operator on $\Ha$
that preserves both the forms $S_1$ and $S_2$.

\bigskip
\noindent
{\bf Proof:}
$M$ is defined by the identity $S_1 (v, M w) = S_2 (v, w)$. If
$T$ preserves boths forms then one has the string of equalities:

$$S_1 (Tv, TM w) =
S_1 (v, M w) =
S_2 (v,  w) =
S_2 (Tv, T w) =
S_1 (Tv, MT w) $$

\bigskip
\nn
Since $T$ is assumed invertible, this is the same as saying

$$ S_1 (v, TM w) =
S_1 (v, MT w) , \ \for \ \all \ v , w \in \Ha
\eqno (56)$$

\bigskip
\nn
But $S_1$ is non-degenerate, namely $\Si_1$ was an
isomorphism. Therefore (56) implies that $TM \equiv MT$, as desired.
\xx

\bigskip
It is clear that to prove $S_2$ is a real multiple of $S_1$
means that the intertwining  operator $M$ has to be simply
multiplication by a scalar. This can now be deduced by looking
at the complexified representation of
M\"ob $\left({  S^1}\right)$ on $\Ha_\CC$, which is unitary,
and applying Schur's Lemma.

\bigskip
\noindent
{\bf LEMMA 4.6:}
The unitary representation of $SL(2, \RR)$ on $\Ha_\CC$
decomposes into precisely two irreducible pieces - namely on
$W_+ $ and $W_-$. In fact these two representations correspond
to the two lowest (conjugate) members in the discrete series
for
$SL (2, \RR)$.

\bigskip
\noindent
{\bf Proof:}
We refer to [10] or [19] for the list of {\it
irreducible unitary representations} of $SL(2, \RR)$ that
constitute what is called its ``discrete series''. Each of
these representations is indexed by an integer $m = \pm 2, \pm
3, \pm 4, \cdots$. For $m \geq 2$ one can write this
representation on the $L^2$ space of holomorphic functions in
$\D$ with the following weighted Poincar\'e measure:

$$d \nu_m = (1 - |z|^2)^m
{dx dy \over (1 - |z|^2)^2 } \ , \ |z| < 1.
\eqno (57)$$

On the Hilbert space $L_{\Hol}^2
(\D, d \nu_m)$ the discrete series representation of
$SL (2, \RR)$ corresponding to this chosen value of $m$ is
given by
$\pi_m:SL(2,\RR) \ra \Aut
\left({ L_{\Hol}^2 (\D, d \nu_m)}\right)$, where

$$\pi_m (\g) (f(z)) = f
\left({ {az + b \over cz + d} }\right)
(cz + d)^{- m}.
\eqno (58)$$

\nn
Here, of course, $\g \in SL(2, \RR)$ corresponds to the
$(PSU(1, 1))$ M\"obius transformation $
{\ds az + b \over \ds cz + d} $ on the disc obained as usual
by conjugating the $SL(2, \RR)$ matrix by the M\"obius
isomorphism (Cayley transform) of the upper half-plane onto the disc.

For $m \leq - 2$ the anti-holomorphic functions conjugate to
those in the above Hilbert spaces need only be used.

We claim that the representation given by the operators
$V_\vp$ on $W_+$ (equation (35)), $\vp \in $ M\"ob
$\left({  S^1}\right)$, can be indentified with the $m = 2$ case
of the above discrete series of representations of $SL(2,
\RR)$. Note, M\"ob
$\left({  S^1}\right) \equiv PSU(1, 1) \cong SL(2, \RR)/(\pm
I)$. Recall from Theorem 2.1 that $W_+$ is identifiable as
$\Hol_2 (\D)$. The action of $\vp$ is given on $\Hol_2$ by :

$$V_\vp (F) = F \circ \vp - F \circ \vp (0), ~~F \in \Hol_2 (\D).
\eqno (59)$$

\nn
But $\Hol_2$ consists of normalized $(F(0) = 0) $ holomorphic
functions in $\D$ whose {\it derivative is in}
$L^2 (\D$, Euclidean measure). From (59), by the chain rule,

$$
{d \over dz } V_\vp (F) =
\left({ { dF \over dz } \circ \vp }\right)  \vp'
\eqno (60)$$

\bigskip
So we can rewrite the representation on the derivatives of the
functions in $\Hol_2$ by the formula (60) - which coincides
with formula (58) for $m = 2$. Indeed $d \nu_2$ is, by (57),
simply the Euclidean (Lebesgue) measure on the disc and thus
$L_{\Hol}^2 (\D, d \nu_2) \cong \Hol_2 (\D)$. (This last isomorphism
being given by sending $F \in \Hol_2 (\D)$ to its derivative.)
Our claim is proved. \xx

\bigskip
It is clear that the representation on the conjugate space will
correspond to the $m = - 2$ (highest weight vector of weight $
- 2$) case of the discrete series. In particular, the
representations we obtain of M\"ob $\left({  S^1}\right)$ by unitary
operators of $W_+ $ and $W_-$ are both {\it irreducible. }

\bigskip
\noindent
{\bf Proof of Theorem 4.3:}
By Lemma 4.5, (the $\CC$-linear extension of ) the
intertwining operator $M$ commutes with every one of the
unitary operators $V_\vp : \Ha_\CC \ra \Ha_\CC$ as $\vp$
varies over M\"ob $\left({  S^1}\right)$.
Since $W_+$ and $W_-$ are the only two invariant subspaces for
all the $V_\vp$, as proved above, it follows that $M$ must map
$W_+$ either to $W_+$ or to $W_-$. Let us first assume the
former case. Then $M$ commutes with all the unitary operators
$V_\vp$ on $W_+$, which we know to be an irreducible
representation. Schur's Lemma says that a unitary
representation will be irreducible if and only if the only
operators that commute with all the operators in
the representation are simply the scalars (see [19,
page 11]). Since $M$ was a real operator to start with, the
scalar must be real.

The alternative assumption that $M$ maps $W_+$ to $W_-$ is
untenable. In fact, if that were so we could replace $M$ by
$M$ followed with complex conjugation. This new $M$ will map
$W_+$ to itself and will again commute with all the $V_\vp$,
hence it must be a scalar. Since the original $M$ arose from a
real operator this scalar can be seen to be real. But scalar
multiplication preserves $W_+ $ - hence the intertwining
operator must map $W_+$ (and $W_-$) to itself.\xx

\bigskip
Our proof is complete. The absolute naturality of the symplectic
form thus established will be utilised in understanding the
$H^{1/2}$ space as a {\it Hilbertian} space, -- namely a space
possessing a fixed
symplectic structure but a large family of compatible complex
structures. See the following sections.

%%% PART 2 (SECTIONS 5 to 10 and REFERENCES ) IN SEPARATE FILE %%%%%%%%
%%% (FULL TEX PAPER TOO BIG FOR E-MAIL MAILER) %%%%%%%%%%%%%%%%%%%%%%%%
%%%%%%%%%%%%%%%%%%%%%%%%%%%%%%%%%%%%%%%%%%%%%%%%%%%%%%%%%%%%%%%%%%%%%%%
%%%%%%%%%%%%%%%%%%%%%%%%%%%%%%%%%%%%%%%%%%%%%%%%%%%%%%%%%%%%%%%%%%%%%%%
%%%%%% TEICHMULLER THEORY AND THE UNIVERSAL PERIOD MAPPING .. II  %%%%%
%%%%%% PART 2 of NAG-SULLIVAN, Sections 5 to 10 and References.   %%%%%
%%%%%%%%%%%%%%%%%%%%%%%%%%%%%%%%%%%%%%%%%%%%%%%%%%%%%%%%%%%%%%%%%%%%%%%

\bigskip
\nn
{\bf \S 5- The $H^{1/2}$ space as first cohomology:}

The Hilbert space $H^{1/2}$, that is the hero of our tale, can be
interpreted as the first cohomology space with real coefficients of the
"universal Riemann surface" -- namely the unit disc -- in a
Hodge-theoretic sense. That will be fundamental for us in explaining the
properties of the period mapping on the universal Teichm\"uller space.

In fact, in the classical theory of the period mapping,
the vector space $H^1(X,\RR)$
plays a basic role, $X$ being a closed orientable topological
surface of genus $g$ to start with. This real vector space comes
equipped with a canonical symplectic structure given by the cup-product
pairing, $S$. Now, whenever $X$ has a complex manifold structure, this
real space $H^1(X,\RR)$ of dimension $2g$ gets endowed with {\it a
complex structure $J$ that is compatible with the cup-pairing $S$}. This
happens as follows: When $X$ is a Riemann surface, the cohomology space
above is precisely the vector space of real harmonic 1-forms on $X$, by
the Hodge theorem. Then the {\it complex structure $J$ is the Hodge
star operator on the harmonic 1-forms}. The compatibility with the cup
form is encoded in the relationships (61) and (62):

$$
S(J\a, J\b)= S(\a, \b),
{}~~~{\rm for~~ all}~~ \a , \b \in  H^1(X,\RR)
\eqno(61)
$$
\noindent
and that, intertwining $S$ and $J$ exactly as in equation (11),
$$
S(\a, J\b) = inner ~~ product(\a, \b)
\eqno(62)
$$
should define a positive definite inner product on $H^1(X,\RR)$.
[In fact, as we will further describe in Section 7, the Siegel disc of
period matrices for genus $g$ is precisely the space of all the
$S$-compatible complex structures $J$.] Consequently, the period mapping
can be thought of as the variation of the Hodge-star complex structure
on the topologically determined symplectic vector space $H^1(X,\RR)$.
See Sections 7 and 8 below.

\bigskip
\nn
{\bf Remark:} Whenever $X$ has a complex structure, one gets an isomorphism
between the real vector space  $H^1(X,\RR)$ and
the $g$ dimensional complex vector space
$H^1(X,\cal O)$, where $\cal O$ denotes the sheaf of germs of holomorphic
functions. That is so because $\RR$ can be considered as a subsheaf of
$\cal O$ and hence there is an induced map on cohomology.
It is interesting to check that this natural map is an
isomorphism, and that the complex structure so induced on
$H^1(X,\RR)$ is the same as that given above by the Hodge star.

\bigskip
For our purposes it therefore becomes relevant to consider, for an {\it
arbitrary} Riemann Surface $X$, {\it the Hodge-theoretic first
cohomology vector space as the space of $L^2$ (square-integrable) real
harmonic 1-forms on $X$}. This real Hilbert space will be denoted
$\Ha(X)$. Once again, in complete generality, this Hilbert space has a
non-degenerate symplectic form $S$ given by the cup (= wedge) product:

$$
S(\p_1, \p_2) =
\int _{}^{}\!
\int_{X}^{}
\p_1 \wedge \p_2
\eqno(63)
$$

\nn
and the Hodge star is the complex structure $J$ of $\Ha(X)$ which is
again compatible with $S$ as per (61) and (62). In fact,
one verifies that the $L^2$ inner
product on $\Ha(X)$ is given by the triality
relationship (62) -- which is the same as (11).

Since in the universal Teichm\"uller theory we deal with the "universal
Riemann surface" -- namely the unit disc $\D$ -- (being the universal
cover of all Riemann surfaces), we require the

\bigskip
\noindent
{\bf PROPOSITION 5.1:} For the disc $\D$, the Hilbert space $\Ha(\D)$ is
isometrically isomorphic to the real Hilbert space $\Ha$ of Section 2.
Under the canonical identification the cup-wedge pairing is the
canonical symplectic form $S$ and the Hodge star becomes the
Hilbert-transform on $\Ha$.

\medskip
\noindent
{\bf Proof:} For every  $\p \in \Ha(\D)$
there exists a unique real harmonic
function $F$ on the disc with $F(0)= 0$  and $dF = \p$. Clearly  then,
$\Ha(\D)$ is isometrically isomorphic to the Dirichlet space $\Da$ of
normalized real harmonic functions having finite energy. But in Section
2 we saw that this space is isometrically isomorphic to $\Ha$ by passing
to the boundary values of $F$ on $S^1$.

If $\p_1 = dF_1$ and $\p_2 = dF_2$ , then integrating $\p_1 \wedge \p_2$
on the disc amounts to, by Stokes' theorem,
$$
\int_{}^{}\!\int_{\D}^{}  dF_1 \wedge dF_2 =
\int_{S^1}^{} F_1 dF_2 = S(F_1, F_2)
$$
as desired.

Finally, let $\p = udx + vdy$ be a $L^2$ harmonic 1-form with $\p = dF$.
Suppose $G$ is the harmonic conjugate of $F$ with $G(0)=0$.
Then $dF + idG$ is a holomorphic 1-form on $\D$ with real part $\p$. It
follows that the Hodge star maps $\p$ to $dG$; hence, under the above
canonical identification of $\Ha(\D)$ with $\Ha$, we see that the
star operator becomes the Hilbert transform. \xx

\noindent
{\bf Remark on the generalised Jacobi variety:}  The complex torus that is
the Jacobi variety of a closed genus $g$ Riemann surface $X$ can be
described as the complex vector space $(H^1(X,\RR), {\rm Hodge~~ star})$
modulo the lattice $H_1(X,\ZZ)$. Indeed, the integral homology group
acts on $H^1(X,\RR)$ as linear functionals by integration of 1-forms on
cycles, and since
$H^1(X,\RR)$  is a Hilbert space, we may canonically identify the dual
space with itself. Thus
$H_1(X,\ZZ)$ appears embedded inside $H^1(X,\RR)$, and the quotient is the
complex torus that is the classical Jacobi variety of the Riemann
surface.

But for the same reasons as above, for an {\it arbitrary} Riemann surface
$X$, $H_1(X,\ZZ)$ does sit inside  the Hodge-theoretic first cohomology
Hilbert space $\Ha(X)$. And this last space carries, as we know, the
Hodge star complex structure. Thus it makes sense to try to define the
generalised Jacobi variety of $X$ as the quotient of this complex
Hilbert space by the "discrete subgroup" $H^1(X,\ZZ)$. For certain
classes of open Riemann surfaces that quotient is a reasonable object,
and we will report on these matters in future articles. For the unit
disc itself then, the generalised Jacobian
{\it is} the Hilbert space
$H^{1/2}$ = $\Ha$ equipped with the Hilbert transform complex structure.

%%%%%%%%%%%%%%%%%%%%%%%%%%%%%%%%%%%%%%%%%%%%%%%%%%%%%%%%%%%%%%%%%%%%%%%
\bigskip
\nn
{\bf \S 6- Quantum calculus and $H^{1/2}$:}

\bigskip
\nn
A.Connes has proposed (see, for example,
[6] and Connes' book "Geometrie Non-Commutatif")
a "quantum calculus" that associates to a function $f$ an operator that should
be considered its quantum derivative -- so that the operator theoretic
properties of this $d^{Q}(f)$ capture the smoothness properties of the
function.  One advantage is that this operator can undergo all the
operations of the functional calculus. The fundamental definition in one
real dimension is

$$
d^{Q}(f)=[J,M_{f}]
\eqno(64)
$$

\nn
where $J$ is the Hilbert transform in one dimension explained in
Section 2, and $M_{f}$ stands for (the generally unbounded) operator
given by multiplication by $f$. One can think of the quantum derivative
as operating (possibly unboundedly) on the Hilbert space $L^{2}(S^1)$ or
on other appropriate function spaces.

\medskip
\nn
{\bf Note:} We
will also allow quantum derivatives to be taken with respect to other
Hilbert-transform like operators; in particular, the Hilbert transform
can be replaced by some conjugate of itself by a suitable automorphism
of the Hilbert space under concern. In that case we will make explicit
the $J$ by writing $d^{Q}_{J}(f)$ for the quantum derivative. See
Section 8 for applications.

\bigskip
As sample results relating the properties of the quantum derivative
with the nature of $f$, we quote:
$d^{Q}(f)$ is a bounded operator on $L^{2}(S^1)$ {\it if and only if}
the function $f$ is of bounded mean oscillation. In fact, the operator
norm of the quantum derivative is equivalent to the BMO norm of $f$.
Again, $d^{Q}(f)$ is a compact operator on $L^{2}(S^1)$ {\it if and only if}
$f$ is in $L^{\Inf}(S^1)$ and has vanishing mean oscillation.
Also, if $f$ is H\"older, (namely in some H\"older class), then
the quantum derivative acts as a compact operator on H\"older.
See [6], [6b].
(Note that the union of all the H\"older classes is both
quasisymmetrically invariant and Hilbert-transform stable. Moreover,
functions that are of bounded variation and H\"older
form a quasisymmetrically invariant subspace of $H^{1/2}$.)
Similarly, the requirement that $f$ is a member of certain Besov spaces
can be encoded in properties of the quantum derivative.

\bigskip
Our Hilbert space $H^{1/2}(\RR)$ has a very simple
interpretation in these terms:

\bigskip
\nn
{\bf PROPOSITION 6.1:} $f \in H^{1/2}(\RR)$  if and only if the operator
$d^{Q}(f)$ is Hilbert-Schmidt on $L^{2}(\RR)$ [or on $H^{1/2}(\RR)$].
The Hilbert-Schmidt norm
of the quantum derivative {\it coincides} with the $H^{1/2}$ norm of
$f$.

\medskip
\nn
{\bf Proof:} Recall that the Hilbert transform on the real line is given
as a singular integral operator with integration kernel $(x-y)^{-1}$.
A formal calculation therefore shows that

$$
(d^{Q}(f))(g)(x) =
\int_{\RR}^{} {{f(x)-f(y)}\over {x-y}} g(y)dy
\eqno(65)
$$

But the above is an integral operator with kernel $K(x,y)=
(f(x)-f(y))/(x-y)$, and such an operator is Hilbert-Schmidt if and only
if the kernel is square-integrable over $\RR^{2}$. Utilising now the
Douglas identity -- equations (24) or (42) -- we are through. \xx

\bigskip
Since the Hilbert transform, $J$, is the standard complex structure on
the $H^{1/2}$ Hilbert space, and since this last space was shown to
allow an action by the quasisymmetric group, $QS(\RR)$, some further
considerations become relevant. Introduce the operator $L$ on 1-forms on
the line to functions on the line by:

$$
(L\vp)(x) =
\int_{\RR}^{} [log \vert x-y \vert]\vp (y) dy
\eqno(66)
$$

One may think
of the Hilbert transform $J$ as operating on
either the space of functions or on the space of 1-forms (by
integrating against the kernel $dx/(x-y)$). Let $d$ as usual denote
total derivative (from functions to 1-forms). Then notice that $L$ above
is an operator that is essentially a smoothing inverse of the
exterior derivative. In fact, $L$ and $d$ are connected to $J$
via the relationships:

$$
d \circ L = J_{1-forms}; ~~~ L \circ d = J_{functions}
\eqno(67)
$$

\bigskip
\nn
{\bf The Quasisymmetrically deformed operators:} Given any
q.s. homeomorphism $h \in QS(\RR)$ we think of it as producing
a q.s. change of structure on the line, and hence we define the
corresponding transformed operators, $L^{h}$ and $J^{h}$ by
$L^{h}= h \circ L \circ h^{-1}$ and $J^{h}= h \circ J \circ h^{-1}$.
($J$ is being considered on functions in $\Ha = H^{1/2}(\RR)$,
as usual.)
The q.s homeomorphism (assumed to be say $C^1$ for the deformation
on $L$), operates standardly on functions and forms by pullback.
Therefore, {\it $J^{h}$ simply stands for the Hilbert transform
conjugated by the symplectomorphism $T_{h}$ of $\Ha$ achieved by
pre-composing by the q.s. homeomorphism $h$.} $J^{h}$ is thus a
new complex structure on $\Ha$.

\smallskip
\nn
{\bf Note:} The complex structures on $\Ha$ of type
$J^{h}$ are exactly those that constitute the image of $T(1)$ by the
universal period mapping. (See Section 8.)
The target manifold, the universal Siegel space,
can be thought of as a space of $S$-compatible complex
structures on $\Ha$.

Let us write the perturbation achieved by $h$ on these operators as the
"quantum brackets":

$$
\{h,L\}=L^{h} - L ;~~~  \{h,J\}=J^{h} - J.
\eqno(68)
$$

Now, for instance, the operator $d \circ \{h,J\}$ is represented
by the kernel $(h \times h)^{*}m - m$ where
$m = dx dy/{(x-y)^{2}}$. For
$h$ suitably smooth this is simply
$d_{y}d_{x}(log [({h(x)-h(y)})/({x-y})])$.
It is well known that
$(h \times h)^{*}m = m$ when $h$ is a M\"obius transformation.
Interestingly, therefore,
on the diagonal ($x=y$), this kernel
becomes ($1/6$ times) the Schwarzian derivative of $h$ (as a
quadratic differential on the line). For the other operators in
the table below the kernel computations are even easier.

Set $K(x,y) = log [({h(x)-h(y)})/({x-y})]$ for convenience.
We have the following table of quantum calculus formulas:

$$
\matrix {{\bf Operator} & {\bf Kernel} & {\bf On~~ diagonal} &
{\bf Cocycle~~on}~~QS(\RR) \cr
\{ h,L\} &  K(x,y) & log (h^\prime) & function-valued \cr
d \circ \{h,L\} & d_{x}K(x,y) & {h'' \over {h'}}dx & 1-form-valued \cr
d \circ \{h,J\} & d_{y}d_{x}K(x,y) & {1 \over 6}Schwarzian(h)dx^{2} &
quadratic-form-valued}
\eqno(69)
$$

The point here is that these operators make sense when $h$ is merely
quasisymmetric. If $h$ happens to be appropriately smooth, we can
restrict the kernels to the diagonal to obtain the respective
nonlinear classical derivatives (affine Schwarzian, Schwarzian, etc.)
as listed in the table above.

\bigskip
\nn
{\bf Remark:} It is worth pointing out that the central
extensions associated  to the three cocycles in the horizontal
lines of the table
above respectively correspond to the subgroups: {\it (i)Translations,
(ii)Affine transformations, {\rm and} (iii)Projective
(M\"obius) transformations}.

%%%%%%%%%%%%%%%%%%%%%%%%%%%%%%%%%%%%%%%%%%%%%%%%%%%%%%%%%%%%%%%%%%%
\bigskip
\nn
{\bf \S 7- The universal period mapping on $T(1)$:}

\bigskip
Having now all the necessary background results behind us, we are
finally  set to move into the theory of the universal period (or
polarisations) map itself.

\bigskip
The Frechet Lie group, $Diff(S^1)$ operating by pullback (=
pre-composition) on smooth functions, had a faithful representation by
bounded symplectic operators on the symplectic vector space $V$ of
equation (1). That induced the natural map $\Pi$ of the homogeneous space
$M=Diff(S^1)/M\ddot{o}b(S^1)$ into Segal's version of the Siegel space of
period matrices.  In [12] [13] we had shown that this map:

$$
\Pi: Diff(S^1)/M\ddot{o}b(S^1) \ra Sp_{0}(V)/U
\eqno(70)
$$

\nn
is {\it equivariant, holomorphic, K\"ahler isometric immersion}, and
moreover that it qualifies as a
{\it generalised period matrix map} (remembering ([15]) that
the domain is a
complex submanifold of the universal space of Riemann surfaces $T(1)$).

{}From the results of Sections 2, 3, and 4, we know that the full
quasisymmetric group, $QS(S^1)$ operates as bounded symplectic operators
on the Hilbert space $\Ha$ that is the completion of the pre-Hilbert
space $V$. The same proof as offered in the articles quoted demonstrates
that the subgroup of $QS$ acting unitarily is the M\"obius subgroup.
Clearly then we have obtained the {\it extension of } $\Pi$
(also called $\Pi$ to save on nomenclature) {\it to the entire
universal Teichm\"uller space}:

$$
\Pi: T(1) \ra Sp(\Ha)/U
\eqno(71)
$$

\bigskip
Let us first exhibit the nature of the complex Banach manifold that is
the target space of the period map (71). This space, which is the
universal Siegel period matrix space, denoted
$\cal S_{\infty}$, has several interesting descriptions:

\medskip
\nn
{\bf (a):}
$\cal S_{\infty}$=
{the space of positive polarizations of the symplectic Hilbert space
$\Ha$ }. Recall ([12], [13], [18]) that a
positive polarization signifies the choice of
a closed complex subspace $W$ in $\Ha_{\CC}$ such that
(i) $\Ha_{\CC} = W \op \oo W$; (ii) $W$ is $S$-isotropic,
namely $S$ vanishes on
arbitrary pairs from $W$; and (iii) $iS(w, \oo w)$ defines the square of
a norm on $w \in W$.

\medskip
\nn
{\bf (b):}
$\cal S_{\infty}$=
{the space of $S$-compatible complex structure operators on $\Ha$ }.
That consists of bounded invertible operators $J$ of $\Ha$ onto itself
whose square is the negative identity and $J$ is compatible with $S$ in
the sense that requirements (61) and (62) are valid.
Alternatively, these are the complex structure operators $J$ on $\Ha$
such that $H(f,g) = S(f,Jg) + iS(f,g)$ is a positive definite Hermitian
form having $S$ as its imaginary part.

\medskip
\nn
{\bf (c):}
$\cal S_{\infty}$=
{the space of bounded operators $Z$ from $W_{+}$ to $W_{-}$ that satisfy
the condition of $S$-symmetry: $S(Z\a,\b)=S(Z\b,\a)$ and are in the unit
disc in the sense  that $(I-Z \oo Z)$ is positive definite}.
The matrix for $Z$ is the "period matrix" of the classical theory.

\medskip
\nn
{\bf (d):}
$\cal S_{\infty}$=
the homogeneous space $Sp(\Ha)/U$; here $Sp(\Ha)$ denotes all bounded
symplectic automorphisms of $\Ha$, and $U$ is the unitary subgroup
defined as those symplectomorphisms that keep the subspace $W_{+}$
(setwise) invariant.

\bigskip
Introduce the {\it Grassmannian} $Gr(W_{+}, \Ha_{\CC})$ of subspaces of type
$W_{+}$ in $\Ha_{\CC}$, which is obviously a complex Banach manifold
modelled on the Banach space of all bounded operators from $W_{+}$ to
$W_{-}$. Clearly,  $\cal S_{\infty}$  is embedded in $Gr$ as a complex
submanifold.  The connections between the above descriptions of the
Siegel universal space are transparent:

\medskip
\nn
(a:b) the positive polarizing subspace $W$ is the $-i$-eigenspace of the
complex structure operator $J$ (extended to $\Ha_{\CC}$ by complex
linearity).

\medskip
\nn
(a:c) the positive polarizing subspace $W$ is the graph of the operator
$Z$.

\medskip
\nn
(a:d) $Sp(\Ha)$ acts transitively on the set of positive polarizing
subspaces. $W_{+}$ is a polarizing subspace, and the isotropy
(stabilizer) subgroup thereat is exactly $U$.

\bigskip
\nn
{\bf $\Ha$ as a Hilbertian space:} Note that the method (b) above
describes the universal Siegel space as a space of Hilbert space
structures on the fixed underlying symplectic vector space $\Ha$. By the
result of Section 4 we know that the symplectic structure on $\Ha$ is
completely canonical, whereas each choice of $J$
above gives a Hilbert space inner
product on the space by intertwining $S$ and $J$ by the fundamental
relationship (11) (or (62)). Thus $\Ha$ is a {\it "Hilbertian space"},
which signifies a complete topological vector space with
a canonical symplectic structure but lots of compatible inner products
turning it into a Hilbert space in many ways.

\bigskip
\nn
We come to one of our {\it Main Theorems:}

\nn
{\bf THEOREM 7.1:} The universal period mapping $\Pi$ is an injective,
equivariant, holomorphic immersion between complex Banach manifolds.

\nn
{\bf Proof:}  From our earlier papers [12] [13] we know these facts
for $\Pi$ restricted to $M$. The proof of equivariance is the same
(and simple). The map is an injection because if we know the subspace
$W_{+}$ pulled back by $w_{\mu}$, then we can recover the q.s.
homeomorphism $w_{\mu}$. In fact, inside the given subspace look at
those functions which map $S^1$ homeomorphically on itself.
One sees easily that these must
be precisely the M\"obius transformations of the circle pre-composed
by $w_{\mu}$. The injectivity (global Torelli theorem) follows.

Let us write down the matrix for the symplectomorphism $T$ on $\Ha_{\CC}$
obtained by pre-composition by $w_{\mu}$. We will write in the standard
orthonormal basis $e^{ik\theta}/k^{1/2}$, $k=1,2,3..$ for $W_{+}$, and
the complex conjugates as o.n. basis for $W_{-}$.

In $\Ha_\CC = W_+ \op W_- $
block form, $T$ is given by maps: $A: W_{+} \ra W_{+}$,
$B:W_{-} \ra W_{+}$. The conjugates of $A$ and $B$ map $W_{-}$ to $W_{-}$
and $W_{+}$ to $W_{-}$, respectively. The matrix entries for
$A=((a_{pq}))$ and $B=((b_{rs}))$ turn out to be:

$$
a_{pq}={(2\pi)^{-1}}{p^{1/2}}{q^{-1/2}}\int_{0}^{2\pi}
{(w_{\mu}(e^{i\theta}))^{q}}{e^{-ip\theta}}d\theta, ~~p,q \ge 1
$$
$$
b_{rs}={(2\pi)^{-1}}{r^{1/2}}{s^{-1/2}}\int_{0}^{2\pi}
{(w_{\mu}(e^{i\theta}))^{-s}}{e^{-ir\theta}}d\theta, ~~r,s \ge 1
$$

Recalling the standard action of symplectomorphisms on the Siegel disc
(model {\bf (c)} above), we see that the corresponding operator
[=period matrix] $Z$ appearing from the Teichm\"uller point [$\mu$]
is given by:

$$
\Pi[\mu] = {\oo B}{A^{-1}}
$$
The usual proof of finite dimensions shows that for any symplectomorphism
$A$ must be invertible -- hence the above explicit formula makes sense.

Since the Fourier coefficients appearing in $A$ and $B$ vary only
{\it real-analytically} with $\mu$, it may be somewhat surprising that
$\Pi$ is actually {\it holomorphic}. In fact,
a computation of the first variation of $\Pi$ at the origin of $T(1)$
( i.e., the derivative map) in the Beltrami direction
$\nu$ shows that the following {\it Rauch variational formula} subsists:

$$
{(d\Pi([\nu]))_{rs}}={{\pi}^{-1}}{(rs)^{1/2}}\int_{}^{}\int_{\Delta}^{}
\nu(z){z^{r+s-2}}dxdy
$$

\nn
The proof of this formula is as shown for $\Pi$ on the smooth
points submanifold $M$ in our earlier papers. The
manifest complex linearity of the derivative, i.e.,
the validity of the Cauchy-Riemann equations, combined with equivariance,
demonstrates that $\Pi$ is complex analytic on $T(1)$, as desired.
\xx

\bigskip
\nn
{\bf Interpretation of $\Pi$ as period map:}
Let us take a moment to recall why the map $\Pi$ qualifies as a
universal version of the classical genus $g$ period maps. As we had
explained in our previous papers, in the light
of P.Griffiths' ideas, the classical period map may be thought of as
associating to a Teichm\"uller point a positive polarizing subspace of
the first cohomology $H^1(X,\RR)$. The point is that when $X$ has a
complex structure, then the complexified first cohomology decomposes as:
$$
H^1(X,\CC)=
H^{1,0}(X) \op H^{0,1}(X)
\eqno(72)
$$
\nn
The period map associates the subspace $H^{1,0}(X)$ -- which is positive
polarizing with respect to the cup-product symplectic form --
to the given complex structure on $X$.  Of course, $H^{1,0}(X)$
represents the holomorphic 1-forms on
$X$, and that is why this is nothing but the usual period mapping.

{\it But that is precisely what $\Pi$ is doing in the universal
Teichm\"uller space.} Indeed, by the results of Section 5, $\Ha$ is the
Hodge-theoretic real first cohomology of the disc, with $S$ being the
cup-product.

The standard complex structure on the unit disc
has holomorphic 1-forms that are of the form $dF$ where $F$ is a
holomorphic function on $\D$ with $F(0)=0$. Thus the boundary values of
$F$ will have only
positive index Fourier modes -- corresponding therefore
to the polarizing subspace
$W_{+}$. Now, an arbitrary point of $T(1)$ is described by the choice of
a Beltrami differential $\mu$ on $\D$ perturbing the complex structure.
We are  now asking for the holomorphic 1-forms on $\D_{\mu}$.
Solving the Beltrami equation on $\D$ provides us with
the $\mu$-conformal quasiconformal self-homeomorphism $w_{\mu}$ of the
disc.  This  $w_{\mu}$  is
a holomorphic uniformising coordinate for  the disc with the $\mu$
complex structure.  The holomorphic 1-forms subspace,
$H^{1,0}(\D_{\mu})$,
should therefore comprise those functions on $S^1$
that are the $W_{+}$ functions {\it precomposed with the boundary
values of the q.c. map $w_{\mu}$.}  That is
exactly the action of $\Pi$ on the Teichm\"uller class of $\mu$.
This explains in some detail why
$\Pi$ behaves as an infinite dimensional period mapping.

\bigskip
\nn
{\bf Remark:} On Segal's $C^{\Inf}$ version of the Siegel space --
constructed using Hilbert-Schmidt operators $Z$, there existed the
universal {\it Siegel symplectic metric}, which we studied in [12] [13] and
showed to be the same as the Kirillov-Kostant (= Weil-Petersson) metric
on $Diff(S^1)/Mob(S^1)$. For the bigger Banach manifold
$\cal S_{\infty}$  above, that pairing fails to converge on arbitrary
pairs of tangent vectors because the relevant operators are not any  more
trace-class in general. The difficulties asociated with this matter will be
addressed in Section 9 below, and in further work that is in progress.

\bigskip
\nn
{\bf \S 8- The universal Schottky locus and quantum calculus:}

\bigskip
\nn
Our object is to exhibit the image of $\Pi$ in
$\cal S_{\Inf}$. The result (equation (73)) can be recognized to be a
quantum "integrability condition" for complex structures on the circle
or the line.

\bigskip
\nn
{\bf PROPOSITION 8.1:}
If a positive polarizing subspace $W$ is in the ''universal Schottky
locus'', namely if $W$ is in the
image of $T(1)$ under the universal period
mapping $\Pi$, then
$W$ possesses a dense subspace which is
{\it multiplication-closed} (i.e., an ``algebra" under pointwise
multiplication modulo subtraction of mean-value.)
In quantum calculus terminology, this means that
$$
[d^{Q}_{J},J] = 0
\eqno(73)
$$
where $J$ denotes the $S$-compatible complex structure of $\Ha$ whose
$-i$-eigenspace is $W$. (Recall the various descriptions of
$\cal S_{\Inf}$ spelled out in the last section.)

\bigskip
\nn
{\bf Multiplication-closed polarizing subspace:}
The notion of being multiplication-closed is well-defined
for the relevant subspaces in $\Ha_{\CC}$.
Let us note that the original polarizing subspace
$W_{+}$ contains the dense subspace of holomorphic trigonometric
polynomials (with mean zero) which constitute an algebra.
Indeed, the identity map of $S^1$ is
a member of $W_{+}$, call it $j$, and positive integral powers of $j$
clearly generate $W_{+}$ -- since polynomials in $j$ form a dense subspace
therein. Now if $W$ is any other positive polarizing subspace, we know
that it is the image of $W_{+}$ under some $T \in Sp(\Ha)$. Thus, $W$
will be multiplication-closed precisely when the image of $j$ by $T$
generates $W$, in the sense that its positive integral powers (minus the
mean values) also lie in
$W$ (and hence span a dense subspace of $W$).

In other words, we are considering $W$ ($\in \cal S_{\Inf}$ [description (a)])
to be multiplication-closed provided that the pointwise products of
functions from $W$ (minus their mean values) that happen to be $H^{1/2}$
functions actually land up in the subspace $W$ again. Multiplying $f$ and
$g$ modulo arbitrary additive constants demonstrates that this notion is
well-defined when applied to a subspace.

\bigskip
\nn
{\bf Quantum calculus and equation (73):} We suggest a quantum version
of complex structures in one real dimension, and note that the
integrable ones correspond to the universal Schottky locus under study.

In the spirit of algebraic geometry one takes the real Hilbert space of
functions $\Ha$ = $H^{1/2}(\RR)$ as the ``coordinate ring" of the real line.
Consequently, a complex structure on $\RR$ will be considered to be a
complex structure on this Hilbert space. Since
$\cal S_{\Inf}$ was a space of (symplectically-compatible) complex
structures on  $\Ha$, we are interpreting
$\cal S_{\Inf}$ as a space of quantum  complex structures on the line
(or circle).

Amongst the points of the universal Siegel space, those that can be
interpreted as the holomorphic function algebra for some complex
structure on the circle qualify as the ``integrable'' ones. But $T(1)$
parametrises all the quasisymmetrically related circles, and for
each one, the map $\Pi$ associates to that structure the holomorphic
function algebra corresponding to it; see the interpretation we provided
for $\Pi$ in the last section. It is clear therefore that $\Pi(T(1)$
should be the integrable complex structures. The point is that taking
the standard circle as having integrable complex structure, all the
other integrable complex structures arise from this one by a $QS$ change
of coordinates on the underlying circle. These are the complex
structures $J^{h}$ introduced in Section 6 on quantum calculus. The
$-i$-eigenspace for $J^{h}$ is interpreted as
the algebra of analytic functions on the
quantum real line with the $h$-structure. We will see in the proof
that (73) encodes just this condition.

\bigskip
\nn
{\bf Proof of Proposition 8.1:}  For a point of
$T(1)$ represented by a q.s. homeomorphism $\p$, the period map sends it
to the polarizing subspace $W_{\p} = W_{+} \circ \p$. But $W_{+}$ was a
multiplication-closed subspace, generated by just the identity map $j$ on
$S^1$, to start with. Clearly then, $\Pi(\p) = W_{\p}$ is also
multiplication-closed in the sense explained, and is generated by the
image of the generator of $W_{+}$ -- namely by the q.s homeomorphism
$\p$ (as a member of $\Ha_{\CC}$).
\xx

\medskip
We suspect that the converse is also true: that the $T(W_{+})$ is such
an ''algebra'' subspace for a symplectomorphism $T$ in $Sp(\Ha)$ only
when $T$ arises as pullback by a quasisymmetric homeomorphism of the circle.
This converse assertion is reminescent of standard theorems in Banach
algebras where one proves, for example, that every (conjugation-preserving)
algebra automorphism of
the algebra $C(X)$ (comprising continuous functions on a
compact Hausdorff space $X$) arises from homeomorphisms of $X$.
[Remark of Ambar Sengupta.]
Owing to the technical hitch that $H^{1/2}$ functions are not in general
everywhere defined on the circle, we are as yet unable to find a
rigorous proof of this converse.

Here is the sketch of an idea for proving the converse.
Thus, suppose we are given a subspace $E$ that is multiplication-closed in the
sense explained. Now, $Sp(\Ha)$ acts transitively on the set of
positive polarizing subspaces. We consider a $T \in Sp(\Ha)$ that maps
$W_{+}$ to $E$ preserving the algebra structure (modulo subtracting off
mean values as usual). Denote by $j$ the identity function on $S^1$
and let $T(j)=w$ be its image in $E$.

Since $j$ is a homeomorphism and $T$ is an invertible real symplectomorphism,
one expects that $w$ is also a homeomorphism on $S^1$.  (Recall the signed
area interpretation of the canonical form (8).) It then follows
that the $T$ is nothing other that precomposition by this
$w$. That is because:
$$
T(j^{m}) = T(j)^{m} - {\rm mean ~ value} =
(w(e^{i\t}))^{m} - {\rm mean ~ value} =
j^{m} \circ w - {\rm mean ~ value.}
$$
\nn
Knowing $T$ to be so on powers of $j$ is sufficient, as polynomials
in $j$ are dense in $W_{+}$.

Again, since $T$ is the complexification of a real symplectomorphism,
seeing the action of $T$ on $W_{+}$ tells us $T$ on all of $\Ha_{\CC}$;
namely, $T$ is everywhere precomposition by that homeomorphism $w$ of
$S^1$. {\it By the necessity part of Theorem 3.1 we see that $w$ must be
quasisymmetric}, and hence that the given subspace $E$ is
the image under $\Pi$ of the Teichm\"uller point determined by $w$
(i.e., the coset of $w$ in $QS(S^1)/M\ddot{o}b(S^1)$).

\bigskip
\nn
{\it Proof of equation (73)}: Let $J$ be {\it any} $S$-compatible complex
structure on $\Ha$, namely $J$ is an arbitrary point of $\cal S_{\Inf}$
(description (b) of Section 7). Let $J_{0}$ denote the Hilbert transform
itself, which is the reference point in the universal Siegel space;
therefore $J = T J_{0} T^{-1}$ for some symplectomorphism $T$ in
$Sp(\Ha)$. The $-i$-eigenspace for $J_{0}$ is, of course, the reference
polarizing subspace $W_{+}$, and the subspace $W$ corresponding to $J$
consists of the functions $(f+i(Jf))$ for all $f$ in $\Ha$.
Now, the pointwise product of two such typical elements of $W$ gives:
$$
(f+i(Jf))(g+i(Jg))=
[fg - (Jf)(Jg)] + i[f(Jg) + g(Jf)]
$$
\nn
In order for $W$ to be multiplication closed the function on the right
hand side must also be of the form $(h + i(Jh))$. Namely, for all relevant
$f$ and $g$ in the real Hilbert space $\Ha$ we must have:

$$
J[fg - (Jf)(Jg)] = [f(Jg) + g(Jf)]
\eqno(74)
$$

Now recall from the concepts introduced in Section 6 that one can
associate to functions $f$ their quantum derivative operators
$d^{Q}_{J}(f)$
which is the commutator of $J$ with the multiplication operator $M_{f}$
defined by $f$. The quantum derivative is being taken with respect to
any Hilbert-transform-like operator $J$ as explained above. But now a
short computation demonstrates that equation (74) is the same as saying
that:
$$
J \circ d^{Q}_{J}(f) = - d^{Q}_{J}(f)
$$
operating by $J$ on both sides shows that this is the same as (73).
That is as desired. \xx

\bigskip
\nn
{\bf Remark:} For the classical period mapping on the Teichm\"uller
spaces $T_g$ there is a way of understanding the Schottky locus in terms
of Jacobian theta functions satisfying the nonlinear K-P equations. In a
subsequent paper we hope to relate the finite dimensional Schottky solution
with the universal solution given above.

\bigskip
\nn
{\bf Remark:} For the extended period-polarizations mapping $\Pi$, the
Rauch variational formula that was exhibited in [12], [13], [13a],
and here in the proof of Theorem 7.1, continues to hold.

%%%%%%%%%%%%%%%%%%%%%%%%%%%%%%%%%%%%%%%%%%%%%%%%%%%%%%%%%%%%%%%%%%%%%
\bigskip
\nn
{\bf \S 9- The Teichm\"uller space of the universal
lamination and Weil-Petersson:}

The Universal Teichm\"uller Space, $T(1)$=$T(\Delta)$,
is a {\it non-separable} complex Banach
manifold that contains, as properly embedded complex submanifolds, all
the Teichm\"uller spaces, $T_g$, of the classical compact
Riemann surfaces of every genus $g$ ($\geq 2$). $T_g$ is $3g-3$
dimensional and appears (in multiple copies) within
$T(\Delta)$ as the Teichm\"uller space $T(G)$ of the
Fuchsian group $G$ whenever $\Delta/G$ is of genus $g$.
The closure of the union of a family of these embedded $T_g$ in
$T(\Delta)$ turns out to be a separable complex submanifold of $T(\Delta)$
(modelled on a separable complex Banach space). That submanifold
can be identified as being itself the Teichm\"uller space of the
"universal hyperbolic lamination" $H_\infty$. We will show that
$T(H_{\infty})$ carries a canonical, genus-independent version of the
Weil-Petersson metric, thus bringing back into play the K\"ahler
structure-preserving aspect of the period mapping theory.

\bigskip
\noindent
{\bf The universal laminated surfaces:}
Let us proceed to explain the nature of the (two possible) "universal
laminations" and the complex structures on these. Starting from any
closed topological surface, $X$,  equipped with a base point,
consider the inverse (directed) system
of all finite sheeted unbranched covering spaces of $X$ by other
closed pointed surfaces. The covering projections are all required to be
base point preserving, and isomorphic covering spaces are identified. The
{\it inverse limit space} of such an inverse system is the "lamination"
-- which is the focus of our interest.

\medskip
\nn
{\it The lamination $E_{\infty}$:}
Thus, if $X$ has genus one, then, of course,
all coverings are also tori, and one
obtains as the inverse limit of the tower a certain compact
topological space -- every path component of which (the laminating
leaves) -- is identifiable with the complex plane. This space $E_{\infty}$
(to be thought of as the "universal Euclidean lamination")
is therefore a fiber space over the original torus $X$ with the fiber
being a Cantor set. The Cantor set corresponds to all the possible
backward strings in the tower with the initial element being the base
point of $X$. The total space is compact since it is a closed subset of
the product of all the compact objects appearing in the tower.

\medskip
\nn
{\it The lamination $H_{\infty}$:}
Starting with an arbitrary $X$ of higher genus clearly produces the
{\it same} inverse limit space, denoted $H_{\infty}$, independent of the
initial genus. That is because given any two surfaces of genus greater
than one, there is always a common covering surface of higher genus.
$H_{\infty}$ is our universal hyperbolic lamination, whose Teichm\"uller
theory we will consider in this section.  For the same
reasons as in the case of $E_{\infty}$, this new lamination is also a
compact topological space fibering over the base surface $X$ with fiber
again a Cantor set. (It is easy to see that
in either case the space of backward
strings starting from any point in $X$ is an uncountable, compact,
perfect, totally-disconnected space -- hence homeomorphic to the Cantor
set.) The fibration restricted to each individual leaf (i.e., path
component of the lamination) is a universal covering projection. Indeed,
notice that the leaves of
$H_{\infty}$ (as well as of $E_{\infty}$) must  all be
simply connected -- since any non-trivial loop
on a surface can be unwrapped in a finite cover. [That corresponds to
the residual finiteness of the fundamental group of a closed surface.]
Indeed, group-theoretically speaking, covering spaces correspond to the
subgroups of the fundamental group. Utilising only normal subgroups
(namely the regular coverings) would give a cofinal inverse system and
therefore the inverse limit would still continue to be the
$H_{\infty}$ lamination. This way of interpreting things allows us to see
that the transverse Cantor-set fiber actually has a group structure.
In fact it is the pro-finite group that is the inverse limit of all the
deck-transformation groups corresponding to these normal coverings.

\bigskip
\noindent
{\bf Complex structures :}
Let us concentrate on the universal hyperbolic lamination
$H_{\infty}$ from now on.
For any complex structure on $X$ there is clearly a complex
structure induced by pullback on each surface of the inverse system, and
therefore $H_{\infty}$ itself inherits a complex structure on each leaf,
so that now biholomorphically each leaf is the Poincare hyperbolic
plane. If we think of a reference complex structure on $X$, then any new
complex structure is recorded by a Beltrami coefficient on $X$, and one
obtains by pullback a complex structure on the inverse limit in the
sense that each leaf now has a complex structure and the Beltrami
coefficients vary continuously from leaf to leaf in the Cantor-set
direction. Indeed, the complex structures obtained in the above fashion
by pulling back to the inverse limit from a complex structure on any
closed surface in the inverse tower, have the special property that the
Beltrami coefficients on the leaves are locally constant in the
transverse (Cantor) direction. These "locally constant" families of
Beltrami coefficients on  $H_{\infty}$ comprise the {\it transversely
locally constant} (written "TLC") complex structures on the lamination.
The generic complex structure on $H_{\infty}$, where all continuously
varying Beltrami coefficients in the Cantor-fiber direction are
admissible, will be a limit of the TLC subfamily of complex structures.

To be precise, a {\it complex structure} on a lamination $L$ is a
covering of $L$ bylamination charts (disc) $\times$
(transversal) so that the overlap homeomorphisms are complex
analytic on the disc direction. Two complex structures are {\it
Teichm\"uller equivalent} whenever they are related to each
other by a homeomorphism that is homotopic to the identity
through leaf-preserving continuous mappings of $L$. For us $L$
is, of course, $H_{\Inf}$. Thus we have defined the set
$T(H_{\Inf})$.

Note that there is a distinguished leaf in our lamination, namely the
path component of the point which is the string of all the base points.
Call this leaf $l$. Note that all leaves are dense in $H_{\infty}$, in
particular $l$ is dense. With respect to the base complex structure the
leaf $l$ gets a canonical identification with
the hyperbolic unit disc $\Delta$.
Hence we have the natural "restriction to $l$"  mapping of the Teichm\"uller
space of $H_{\infty}$ into the Universal Teichm\"uller space $T(l) = T(1)$.
Since the leaf is dense, the complex structure on it records the entire
complex structure of the lamination.
The above restriction map is therefore actually
{\it injective} (see [20]) and therefore describes
$T(H_{\infty})$ as an embedded complex analytic submanifold in $T(1)$.

Indeed, as we will explain in detail below, $T(H_{\infty})$ embeds as
precisely the closure in $T(1)$ of the union
of the Teichm\"uller spaces $T(G)$ as
$G$ varies over all finite-index subgroups of a fixed cocompact Fuchsian
group. These
finite dimensional classical Teichm\"uller spaces lying within the
separable,  infinite-dimensional
$T(H_{\infty})$, comprise the TLC points of $T(H_{\infty})$.

Alternatively, one may understand the set-up at hand by looking at the
direct system of maps between Teichm\"uller spaces that is obviously
induced by our inverse system of covering maps. Indeed, each covering
map provides an immersion of the Teichm\"uller space of the covered
surface into the Teichm\"uller space
of the covering surface induced by the
standard pullback of complex structure. These immersions are
Teichm\"uller metric preserving, and provide a direct system whose
direct limit, when completed in the Teichm\"uller metric, gives
produces again $T(H_{\infty})$. The direct limit already contains the
classical Teichm\"uller spaces of closed Riemann surfaces, and the
completion corresponds to taking the closure in $T(1)$.

Let us elaborate somewhat more on these various possible
embeddings of $T(H_{\infty})$ [ which is to be thought of as the {\it
universal Teichm\"uller space of {\bf compact} Riemann surfaces}]
within the classical universal Teichm\"uller space $T(\Delta)$.

\medskip
\noindent
{\bf Explicit realizations of $T(H_\infty)$ within the universal
Teichm\"uller space:}
Start with any cocompact (say torsion-free) Fuchsian group $G$
operating on the unit disc $\Delta$,
such that the quotient is a Riemann surface $X$ of {\it arbitrary} genus $g$
greater than one.
Considering the inverse limit of the directed system of all unbranched
finite-sheeted pointed covering spaces over $X$ gives us a copy of the
universal laminated space $H_\infty$ equipped with a complex structure induced
from that on $X$. Every such choice of $G$ allows us to embed the
separable Teichm\"uller space
$T(H_\infty)$
holomorphically in the Bers universal Teichm\"uller space $T(\Delta)$.

To fix ideas, let us think of the universal Teichm\"uller space as:
$ T(\Delta) = T(1) = QS(S^1)/Mob(S^1)$ as usual.

For any Fuchsian group $\Gamma$ define:
$$
QS(\Gamma) = \{ w \in QS(S^1): w\Gamma w^{-1}~is~again~a~Mobius~group.\}
$$
We say that the quasisymmetric homeomorphisms in $QS(\Gamma)$ are
those that are {\it compatible} with $\Gamma$.
Then the Teichm\"uller space $T(\Gamma)$ is $QS(\Gamma)/Mob(S^1)$ clearly
sits embedded within $T(1)$. [We always think of points of $T(1)$ as
left-cosets of the form $Mob(S^1)\circ w$ = $[w]$ for arbitrary
quasisymmetric homeomorphism $w$ of the circle.]

Having fixed the cocompact Fuchsian group $G$, the Teichm\"uller space
$T(H_\infty)$
is now the closure in $T(1)$ of the direct limit of all the Teichm\"uller
spaces $T(H)$ {\it as $H$ runs over all the finite-index subgroups of the
initial cocompact Fuchsian group $G$}. Since each $T(H)$ is actually
embedded injectively within the universal Teichm\"uller space, and since
the connecting maps in the directed system are all inclusion maps, we see
that the direct limit (which is in general a quotient of the disjoint
union) in this situation is nothing other that just the set-theoretic
{\it union of all the embedded $T(H)$ as $H$ varies over all finite index
subgroups of $G$ }. This
union in $T(1)$ constitutes the dense ``TLC" (transversely locally constant)
subset of $T(H_\infty)$. Therefore, {\it the TLC subset of this embedded
copy of $T(H_\infty)$ comprises the M\"ob-classes of all those
QS-homeomorphisms that are compatible with {\it some} finite index
subgroup in $G$.}

We may call the above realization of
$T(H_\infty)$ as  {\it `` the $G$-tagged embedding" of
$T(H_\infty)$} in $T(1)$.

Remark: We see above, that just as the Teichm\"uller space of Riemann
surfaces of any genus $p$ have lots of realizations within the universal
Teichm\"uller space (corresponding to choices of reference cocompact
Fuchsian groups of genus $p$), the Teichm\"uller space of the lamination
$H_\infty$ also has many different realizations within $T(1)$.

Therefore, in the Bers embedding of $T(1)$, this realization
of $T(H_\infty)$ is the intersection of the domain $T(1)$ in the
Bers-Nehari Banach space $B(1)$ with the separable Banach subspace that
is the inductive (direct) limit of the subspaces $B(H)$ as $H$ varies
over all finite index subgroups of the Fuchsian group $G$.
(The inductive limt topology will give a complete (Banach) space; see,
e.g., Bourbaki's "Topological Vector Spaces".) It is relevant to recall that
$B(H)$ comprises the bounded holomorphic quadratic forms for the group
$H$. By Tukia's results, the Teichm\"uller space of $H$ is exactly the
intersection of the universal Teichm\"uller space with $B(H)$.

\nn
{\bf Remark:} Indeed one expects that the
various $G$-tagged embeddings of $T(H_{\Inf}$
must be sitting in general discretely separated from each other in
the Universal Teichm\"uller space. There is a result to this effect
for the various copies of $T(\Gamma)$, as the base group is varied, due to
K.Matsuzaki (preprint -- to appear in Annales Acad. Scient. Fennicae).
That should imply a similar discreteness for the family of embeddings
of $T(H_{\Inf})$ in $T(\Delta)$.

It is not hard now to see how many different copies of the Teichm\"uller
space of genus $p$ Riemann surfaces appear embedded within the
$G$-tagged embedding of $T(H_\infty)$. That corresponds to non-conjugate
(in $G$) subgroups of $G$ that are of index $(p-1)/(g-1)$ in $G$. This
last is a purely topological question regarding the fundamental group of
genus $g$ surfaces.

Modular group: One may look at those elements of the full universal
modular group $Mod(1)$ [quasisymmetric homeomorphism acting by right
translation (i.e., pre-composition) on $T(1)$] that preserve setwise the
$G$-tagged embedding of $T(H_\infty)$. Since the modular group
$Mod(\Gamma)$ on $T(\Gamma)$ is induced by right translations by those
QS-homeomorphisms that are in the normaliser of $\Gamma$:
$$
N_{qs}(\Gamma)=\{t \in QS(\Gamma): t\Gamma t^{-1}=\Gamma \}
$$
it is not hard to see that only the elements of $Mod(G)$ itself will
manage to preserve the $G$-tagged embedding of $T(H_\infty)$.
[Query: Can one envisage some limit of the modular groups
of the embedded Teichm\"uller spaces as acting on $T(H_\infty)$?]

\bigskip
\noindent
{\bf The Weil-Petersson pairing:}
In [20], it has been shown that the tangent (and the cotangent)
space at any point of $T(H_{\Inf})$
consist of certain holomorphic quadratic
differentials on the universal lamination $H_{\Inf}$.
In fact, the Banach space $B(c)$ of tangent holomorphic quadratic
differentials at the Teichm\"uller point represented by the
complex structure $c$ on the lamination, consists
of holomorphic quadratic differentials on the leaves that vary
continuously in the transverse Cantor-fiber direction. Thus
locally, in a chart, these objects look like
$\vp(z,\lambda)dz^{2}$ in self-evident notation; ($\lambda$
represents the fiber coordinate). The lamination $H_{\Inf}$ also
comes equipped with an invariant transverse measure on the
Cantor-fibers (invariant with respect to the holonomy action of
following the leaves). Call that measure (fixed up to a scale)
$d\lambda$. [That measure appears as the limit of (normalized) measures
on the fibers above the base point that assign (at each finite
Galois covering stage) uniform weights to the points in the fiber.]
{}From [20] we have directly therefore our present goal:

\medskip
\nn
{\bf PROPOSITION 9.1:} The Teichm\"uller space $T(H_{\Inf})$ is a
separable complex Banach manifold in $T(1)$ containing the
direct limit of the classical Teichm\"uller spaces as a dense
subset. The Weil-Petersson metrics on the classical $T_{g}$,
normalized by a factor depending on the genus, fit together and
extend to a finite Weil-Petersson inner product on
$T(H_{\Inf})$ that is defined by the formula:

$$
\int_{H_{\Inf}}^{}
\vp_{1} \vp_{2}(Poin)^{-2}dz \wedge d \oo z d\lambda
\eqno(75)
$$
where $(Poin)$ denotes the Poincare conformal factor for the Poincare
metric on the leaves (appearing as usual for all Weil-Petersson formulas).
\xx

\medskip
\nn
{\bf Remark on Mostow rigidity for $T(H_{\Inf})$:} The quasisymmetric
homeomorphism classes comprising this Teichm\"uller space are again very
non-smooth, since they appear as limits of the fractal q.s. boundary
homeomorphisms corresponding to deformations of co-compact Fuchsian
groups. Thus, the transversality proved in [15, Part II] of the finite
dimensional Teichm\"uller spaces with the coadjoint orbit homogeneous
space $M$ continues to hold for $T(H_{\Inf})$. As explained there, that
transversality is a form of the Mostow rigidity phenomenon.
The formal Weil-Petersson
converged on $M$ and coincided with the Kirillov-Kostant metric, but
that formal metric fails to give a finite pairing on the tangent spaces
to the finite dimensional $T_{g}$. Hence the interest in the above
Proposition.

%%%%%%%%%%%%%%%%%%%%%%%%%%%%%%%%%%%%%%%%%%%%%%%%%%%%%%%%%%%%%%%%%%
\bigskip
\nn
{\bf \S 10-The Universal Period mapping and the Krichever map:}

We make some remarks on the relationship of $\Pi$ with the Krichever mapping
on a certain family of Krichever data. This could be useful in developing
infinite-dimensional theta functions that go hand-in-hand with our infinite
dimensional period matrices.

The positive polarizing subspace, $T_{\mu}(W_{+})$, that is assigned by the
period mapping $\Pi$ to a point $[\mu]$ of the universal Teichm\"uller space
has a close relationship with the Krichever subspace of $L^{2}(S^1)$ that
is determined by the Krichever map on certain Krichever data, when $[\mu]$
varies in (say) the Teichm\"uller space of a compact Riemann surface with one
puncture (distinguished point). I am grateful to Robert Penner for
discussing this matter with me.

Recall that in the Krichever mapping one takes a compact Riemann surface
$X$, a point $p \in X$, and a local holomorphic coordinate around $p$ to
start with (i.e., a member of the ``dressed moduli space''). One also
chooses a holomorphic line bundle $L$ over $X$
and a particular trivialization of $L$ over the
given ($z$) coordinate patch around $p$. We assume that the $z$ coordinate
contains the closed unit disc in the $z$-plane. To such data, the Krichever
mapping associates the subspace of $L^{2} (S^1)$ [here $S^1$ is the
unit circle in the $z$ coordinate]
comprising functions which are restrictions to that
circle of holomorphic sections of $L$ over the punctured surface
$X-\{p\}$.

If we select to work in a Teichm\"uller space $T(g,1)$ of pointed
Riemann surfaces of genus $g$, then one may choose $z$ canonically as a
certain horocyclic coordinate around the point $p$. Fix $L$ to be the canonical
line bundle $T^{*}(X)$ over $X$ (the compact Riemann surface). This has
a corresponding trivialization via ``$dz$''. The Krichever image of this
data can be considered as a subspace living on the unit horocycle
around $p$. That horocycle can be mapped over to the boundary circle
of the universal covering disc for $X-\{p\}$ by mapping out by the
natural pencil of Poincare geodesics having one endpoint at a parabolic
cusp corresponding to $p$.

{\it We may now see how to recover the Krichever subspace (for this restricted
domain of Krichever data) from the subspace in
$H^{1/2}_{\CC}(S^1)$ associated to $(X,p)$ by $\Pi$.} Recall that the
functions appearing in the $\Pi$ subspace are the boundary values of
the Dirichlet-finite harmonic functions whose derivatives give the
holomorphic Abelian differentials of the Riemann surface. Hence, to get
Krichever from $\Pi$ one takes Poisson integrals of the functions in the
$\Pi$ image, then takes their total derivative in the universal covering
disc, and restricts these to the horocycle around $p$ that is
sitting inside the universal cover
(as a circle tangent to the boundary circle of the Poincare disc).

Since Krichever data allows one to create the tau-functions of the $KP$
-hierarchy by the well-known theory of the Sato school (and the Russian
school), one may now use the tau-function from the Krichever data to
associate a tau (or theta) function to such points of our universal
Schottky locus. The search for natural theta functions associated to
points of the universal Siegel space $\cal S$, and their possible use
in clarifying the relationship between the universal and classical Schottky
problems, is a matter of interest that we are pursuing.

%%%%%%%%%%%%%%%%%%%%%%%%%%%%%%%%%%%%%%%%%%%%%%%%%%%%%%%%%%%%%%%%%%%%

\bigskip
\centerline {{\bf REFERENCES} }

%%%%%%%%%%%%%%%%%%%%%%%%%%%%%%%%%%%%%%%%%%%%%%%%%%%%%%%%%%%%%%%%%%%%
\medskip
\item{1}
L.Ahlfors, Remarks on the Neumann-Poincare integral equation,
{\it Pacific J.Math}, 2(1952),  271-280.

\medskip
\item{2}
L.Ahlfors, {\it Conformal Invariants: Topics in Geometric Function
Theory}, Mc Graw Hill,

\medskip
\item{3}
L. Ahlfors, {\it Lectures on quasiconformal mappings,}
Van Nostrand, (1966).

\item{4}
A. Beurling and L. Ahlfors,
The boundary correspondence for quasiconformal mappings,
{\it Acta Math}. 96 (1956), 125-142.

\medskip
\item{5}
A. Beurling and J.Deny, Dirichlet Spaces, {\it Proc. Nat. Acad. Sci.}, 45
(1959), 208-215.

\medskip
\item{6}
A. Connes and D. Sullivan, Quantum calculus on $S^1$ and Teichm\"uller
theory, IHES preprint, 1993.

\medskip
\item{6a}
A.Connes, {\it Geometrie Non Commutative}.

\medskip
\item{6b}
R.R. Coifman, R. Rochberg and G. Weiss, Factorization theorems for Hardy
spaces in several variables, {\it Annals Math.} 103 (1976),611-635.

\medskip
\item{7}
J. Deny and J.L. Lions,
Les espaces de type de Beppo-Levi,
{\it Annales Inst. Fourier}, 5 (1953-54), 305-370.

\medskip
\item{8}
F. Gardiner and D. Sullivan,
Symmetric structures on a closed curve, {\it American J. Math.}
114 (1992), 683-736.

\medskip
\item{9}
P. Griffiths, Periods of integrals on algebraic manifolds, {\it
Bull. American Math. Soc.,}, 75 (1970), 228-296.

\medskip
\item{10}
S. Lang, $SL_2 (\RR)$, Springer-Verlag, (1975).

\medskip
\item{11}
O. Lehto and K. Virtanen,
{\it Quasiconformal mappings in the plane,}
2nd ed. Springer-Verlag, Berlin, (1973).

\medskip
\item{12}
S. Nag,
A period mapping in universal Teichm\"uller space,
{\it Bull. American Math. Soc.,} 26, (1992), 280-287.

\medskip
\item{13}
S. Nag,
Non-perturbative string theory and the diffeomorphism group of
the circle, in {\it ``Topological and Geometrical Methods in
Field Theory''} Turku, Finland International Symposium,
(eds. J. Mickelsson and O.Pekonen), World Scientific, (1992).

\medskip
\item{13a}
S. Nag,
On the tangent space to the universal Teichm\"uller space, {\it Annales
Acad. Scient. Fennicae}, vol 18, (1993), (in press).

\medskip
\item{14}
S. Nag,
{\it The complex analytic theory of Teichm\"uller spaces,}
Wiley-Interscience, (1988).

\medskip
\item{15}
S. Nag and A. Verjovsky,
Diff $(S^1)$ and the Teichm\"uller spaces,
{\it Commun. Math. Physics}, 130 (1990), 123-138 (Part I by
S.N. and A.V. ; Part II by S.N.).

\medskip
\item{16}
M.S. Narasimhan,
The type and the Green's kernel of an open Riemann surfaces,
{\it Annales Inst. Fourier,} 10 (1960), 285-296.

\medskip
\item{17}
H. Reimann, Ordinary differential equations and quasiconformal mappings,
{\it Inventiones Math.} 33 (1976), 247-270.

\medskip
\item{18}
G. Segal,
Unitary representations of some infinite dimensional groups,
{\it Commun. in Math. Physics} 80 (1981), 301-342.

\medskip
\item{19}
M. Sugiura,
{\it Unitary representations and harmonic analysis - an
introduction, } 2nd ed., North Holland/Kodansha, (1975)

\medskip
\item{20}
D. Sullivan, Relating the universalities of Milnor-Thurston, Feigenbaum
and Ahlfors-Bers, Milnor Festschrift volume ``Topological methods
in modern Mathematics'', (ed. L.Goldberg and A. Phillips),
Publish or Perish, 1993, 543-563.

\medskip
\item{21}
H. Triebel,
{\it Theory of function spaces,}
Birkh\"auser-Verlag, (1983).

\medskip
\item{22}
E. Witten,
Coadjoint orbits of the Virasoro group,
{\it Commun. in Math. Physics,} 114 (1981), 1-53.

\medskip
\item{23}
A. Zygmund,
{\it Trigonometric Series,}
Vol 1 and 2,
Cambridge Univ. Press (1968).

%%%%%%%%%%%%%%%%%%%%%%%%%%%%%%%%%%%%%%%%%%%%%%%%%%%%%%%%%%%%%%%%%%%%
\bigskip
\centerline{\bf Authors' addresses}

\bigskip

$$\vbox{
\offinterlineskip
\halign{
#               &#   \cr
\tvi\cc{ The Institute of Mathematical Sciences}
& \cc {City Univ. of New York, Einstein Chair}
\cr
\tvi \cc {Madras   600 113, INDIA}
&\cc {New York, N.Y. 10036, U.S.A.}
\cr
\tvi \cc {and}
&
\cc {and}
\cr
\tvi \cc{I.H.E.S.}
&
\cc{I.H.E.S.}\cr
\tvi \cc {91440 Bures sur Yvette, FRANCE}
&
\cc{91440 Bures sur Yvette, FRANCE}
\cr
}}$$

\end